\documentclass[a4paper]{amsart}

\usepackage{color}
\usepackage{graphicx,epstopdf}
\usepackage{amsmath}
\usepackage{amssymb}
\usepackage{amsfonts}
\usepackage{amsthm}

\theoremstyle{remark}
\newtheorem{remark}{Remark}[section]
 
\begin{document}
 \title[]{Multi-domain spectral approach for the Hilbert transform on 
 the real line}

\author{Christian Klein$^{*}$}
\address{Institut de Math\'ematiques de Bourgogne, UMR 5584\\
                Universit\'e de Bourgogne-Franche-Comt\'e, 9 avenue Alain Savary, 21078 Dijon
                Cedex, France\\
    E-mail Christian.Klein@u-bourgogne.fr}
\author{Julien Riton}
\address{Institut de Math\'ematiques de Bourgogne, UMR 5584\\
                Universit\'e de Bourgogne-Franche-Comt\'e, 9 avenue Alain Savary, 21078 Dijon
                Cedex, France\\
    E-mail julien.riton@gmail.com}

\author{Nikola Stoilov}
\address{Institut de Math\'ematiques de Bourgogne, UMR 5584\\
                Universit\'e de Bourgogne-Franche-Comt\'e, 9 avenue Alain Savary, 21078 Dijon
                Cedex, France\\
    E-mail Nikola.Stoilov@u-bourgogne.fr}
\date{\today}

\begin{abstract}
    A  multi-domain spectral method is presented to compute the Hilbert 
    transform on the whole compactified real line, with a special focus on piece-wise 
	analytic functions and functions with algebraic decay towards infinity. Several examples of these and other types of functions are 
    discussed. As an application solitons to generalized Benjamin-Ono 
	equations are constructed.  
\end{abstract}

 
\thanks{This work was partially supported by 
the ANR-FWF project ANuI - ANR-17-CE40-0035, the isite BFC project 
NAANoD, the ANR-17-EURE-0002 EIPHI and by the 
European Union Horizon 2020 research and innovation program under the 
Marie Sklodowska-Curie RISE 2017 grant agreement no. 778010 IPaDEGAN. 
We thank 
J.A.C.~Weideman for helpful discussions and hints}
\maketitle

\section{Introduction}

We present an efficient numerical approach based on a multi-domain 
spectral method for the computation of the Hilbert 
transform on the real line. We are specifically interested in 
functions which are piece-wise analytic on $\mathbb{R}$, but we also discuss various other examples. The Hilbert transform of a function $f\in L^{2}(\mathbb{R})$ 
 is defined as 
\begin{equation}
    \mathcal{H}[f](x):=\frac{1}{\pi}\mathcal{P}\int_{\mathbb{R}}^{}\frac{f(y)}{x-y}dy
    \label{hilbert},
\end{equation}
where $\mathcal{P}$ denotes the principal value.  The Hilbert transform appears in countless applications in mathematics, physics and signal processing. 
Some important examples include singular integral equations, see 
e.g.~\cite{Mus} where the Hilbert transform  is used as the Cauchy integral on the real line. It is fundamental in linear 
response theory in the form of the 
Kramers-Kronig relations, for applications see \cite{HH09}. Our main interest is in theory of water 
waves where the Hilbert transform appears for instance in the context 
of the generalized
Benjamin-Ono (BO)
equation,
\begin{equation}
    u_{t} +u^{m-1}u_{x} -\mathcal{H}u_{xx} = 0, 
    \label{BO}
\end{equation}
where $m=2,3,\ldots$, see \cite{sautBO} for a recent review and 
\cite{RWY} for a numerical study. 

A convenient way to compute the Hilbert transform is via its 
Fourier transform, defined for
a function $f\in L^{2}(\mathbb{R})$  as
\begin{equation}
    \hat{f}(k) = \mathcal{F}f := \int_{\mathbb{R}}^{}f(y)e^{-iky}dy
    \label{fourier},
\end{equation}
where $k\in \mathbb{R}$ is the dual variable to $x$. It is well known 
that the Fourier symbol of $\mathcal{H}$ 
is simply given by 
\begin{equation}
    \mathcal{F}\mathcal{H}=-i\text{ sgn}(k)
    \label{fourierhilbert},
\end{equation}
i.e., it is not smooth. With a Paley-Wiener type argument this 
immediately implies that the Hilbert transform $\mathcal{H}(f)$ cannot be rapidly 
decreasing in $x$  for $|x|\to\infty$ even for functions \( f \) in the 
Schwartz class $\mathcal{S}(\mathbb{R})$ of rapidly decreasing smooth functions because otherwise its Fourier transform would be smooth.

A standard numerical approach to compute the Hilbert transform is based on 
an approximation of the Fourier transform by a discrete Fourier 
transform (DFT). This is a spectral method, i.e., the numerical error in 
approximating analytic periodic functions decreases exponentially 
with the number $N_{\mathcal{F}}$ of Fourier modes. In addition, the discrete 
Fourier transform can be efficiently computed with the fast Fourier 
transform (FFT) which is known to take $\mathcal{O}(N_{\mathcal{F}}\ln N_{\mathcal{F}})$ operations instead of the $\mathcal{O}(N_{\mathcal{F}}^2)$ the direct implementation of DFT takes. 
Thus, for functions in the Schwartz class, which can be 
seen as smooth and periodic on sufficiently large periods 
within the finite numerical precision, such methods are highly efficient. The 
problem in the context of the Hilbert transform is the singular symbol 
(\ref{fourierhilbert}) which, as mentioned above, implies that the 
Hilbert transform decreases only algebraically in $1/|x|$ for 
$|x|\to\infty$. Such functions are not efficiently approximated via Fourier series.

Weideman \cite{weideman} gave an elegant way to overcome these 
problems by introducing a mapping of the whole real line to the 
circle.
This allows to take advantage of the efficiency of the FFT 
whilst avoiding the disadvantages of the approximation of Hilbert 
transforms via trigonometric functions in the integration variable \( y \). The method, together with rigorous error analysis, is illustrated for 
several examples in \cite{weideman}. There are many other numerical approaches for the 
computation of the Hilbert transform, see for instance
\cite{King} for a recent review and \cite{BP,BVM,Pot,ZYLY} for new 
developments. Some of these approaches compute the 
Hilbert transform in terms of certain transcendental functions which 
then have to be computed as well. As we will show in this paper, for 
piece-wise analytic functions it is possible to compute the Hilbert 
transform in terms of elementary polynomials to the order of machine 
precision\footnote{Note, however, that the presented algorithm requires only 
piece-wise differentiable functions, it is just most efficient for 
piece-wise smooth functions}.

In this paper we address potential problems related to the mapping of the whole 
compactified real line to the sphere or a single interval. First, if 
the function of interest $f$ is only analytic on a finite number of 
intervals $I_{n}$, $n=1,..,M$, where $\cup_{n=1}^{M} 
I_{n}=\mathbb{R}\cup\{\infty\}$ and $f\in C^{r}(\mathbb{R})$, $r\geq 
0$, a 
spectral approach will be only of finite order  on the whole 
real line, but of spectral accuracy on each of the intervals $I_{n}$, 
$n=1,\ldots,M$. As an example of such a function Weideman \cite{weideman} considers $f(y)=\exp(-|y|)$ , which we will also discuss 
here. This function is also considered in \cite{olver} by 
mapping the half-lines $\mathbb{R}^{\pm}$ separately. Second, a 
multi-domain method offers the possibility to allocate resolution 
where needed. For instance for rapidly decreasing functions, no 
collocation points will be needed where the functions vanish with 
numerical precision. 

Our multi-domain approach consists thus in mapping each of the intervals $I_{n}$, 
$n=1,\ldots,M$ to the interval $[-1,1]$. For infinite intervals we 
use $1/y$ as a local coordinate near infinity. The integrals are then 
computed with a spectral quadrature scheme as Clenshaw-Curtis 
\cite{CC}. For $N$  collocation points in total the computational cost for the Clenshaw-Curtis quadrature 
is of the order $\mathcal{O}(N^{2}/M)$ and thus of higher complexity 
than 
Weideman's FFT based algorithm, for the same  total number $N_{\mathcal{F}} = N$ of collocation 
points used for his approach. However, we will show that for some of 
the examples of \cite{weideman}, a quadrature based approach is competitive as the total number $N$ of collocation points on 
all of the intervals $I_{n}$, $n=1,\ldots,M$ can be chosen much 
smaller than the $N_{\mathcal{F}}$ necessary for  in the Fourier approach and still achieve machine precision. 

The paper is organized as follows: in Section~2 we introduce our 
approach for functions which have an algebraic decrease for 
$|y|\to\infty$. For the sake of simplicity we discuss in 
detail the case of two intervals. In Section 3 we discuss functions 
which are piece-wise analytic. In section 4 we address
the case of functions with essential singularities at infinity. As an 
application of the approach, we construct solitary waves for 
generalized BO equations (\ref{BO}) in section 5. We add some 
concluding remarks in section 6. 

\noindent\textbf{Notation:}\\
For clarity we establish the following convention: $x$ is the 
external variable for the Hilbert transform, $y$ is the internal, 
thus $f$ generally is defined over $y$ and $\mathcal{H}[f]$ is a 
function of $x$;  $k$ is the standard dual variable of the Fourier transform. The spaces in which $x$ and $y$ live normally coincide so if we e.g. integrate by parts we write $f(x)$ without further discussion.

\section{Hilbert transform for functions analytical on the 
compactified real line}
In this section, we consider functions with an algebraic decay for 
$|y|\to\infty$. The approach is set up for $N_M$ domains, in 
general $N_{M-1}$ finite ones and one infinite one. The choice of the 
number of domains is imposed by the problem we are studying. 
This means that if the function appearing in the Hilbert 
transform is piece-wise smooth (or at least $C^{1}$ for our algorithm), 
say on $N_{M-1}$ finite intervals, these 
intervals are a natural choice for the intervals in the method. In 
addition it can be that the conditioning of the integration scheme, 
for instance Clenshaw-Curtis, which is of the order of 
$\mathcal{O}(N^{2}_{n})$, where $N_{n}$ is the number of collocation 
points in the interval $I_{n}$, see 
\cite{trefethen}, becomes important if $N_{n}$ has to be chosen large. 
As we will later discuss in more detail, see also \cite{CC}, the choice of 
$N_{n}$ depends on the highest Chebyshev  coefficients since they are an 
indicator of the numerical accuracy 
(see the discussion of the examples). This means that the Chebyshev  
coefficients $c_{n}$, see (\ref{chebcoeff}),  on each considered interval should be of the 
wanted order, say $10^{-16}$, for $n\sim N_{n}$ on each interval $I_{n}$. In  
cases where the number of collocation points in the $n$-th interval 
$N_n$ would have to be chosen too large (in practice  much larger than $100$),  it can be 
beneficial to subdivide this interval  into several intervals such 
that each new $N_n$ can be 
chosen small (in practice around 100). For the 
ease of presentation, we discuss in detail below the case with two intervals one of  which is infinite. 

\subsection{Finite intervals}
We first address the case of a finite interval 
$I_{n}=[a_{n},b_{n}]$, $\ldots<a_{n}<b_{n}<a_{n+1}<\ldots$, $n\leq 
N_{M-1}$.  This means we consider the integral 
\begin{equation}
H_n(x) =  \mathcal{P}\int_{a_{n}}^{b_{n}}\frac{f(y)}{x-y}dy
    \label{I1}.
\end{equation}
If $x\notin I_{n}$, the integrand of (\ref{I1})  is regular, and 
standard quadrature formulae could be applied directly. If \( x\in 
I_{n} \), 
the principal value for $x\in I_{n}$ can be computed in classical manner,
\begin{align}
    H_n(x) = \mathcal{P}\int_{a_{n}}^{b_{n}}\frac{f(y)}{x-y}dy& = 
    \int_{a_{n}}^{b_{n}}\frac{f(y)-f(x)}{x-y}dy
     +f(x)\mathcal{P}\int_{a_{n}}^{b_{n}}\frac{1}{x-y}dy
    \nonumber\\
     & =\int_{a_{n}}^{b_{n}}\frac{f(y)-f(x)}{x-y}dy
     -f(x)\ln\frac{b_{n}-x}{x-a_{n}}
    \label{I2}.
\end{align}
 
Note that 
    the appearance of a logarithm in (\ref{I2}) does not imply that 
    the Hilbert transform is unbounded since there will be similar 
    terms from the other intervals $I_{n}$ leading to a possibly regular 
    expression on the whole real line (depending on the regularity of 
	\( f \), see the example in subsection~ 
    \ref{infinite}).

To compute the regular integrals in (\ref{I1}) or (\ref{I2}), we map 
them to the interval $[-1,1]$ 
via  $y=b_{n}(1+l)/2+a_{n}(1-l)/2$, where $l\in[-1,1]$. 
The integrals of the form 
\begin{equation}
    \int_{-1}^{1}g(l)dl
    \label{gs}
\end{equation}
are computed with the Clenshaw-Curtis algorithm \cite{CC}: we 
introduce the standard Chebyshev collocation points 
\begin{equation}
    l_{m} = \cos(m\pi/N_{n}),\quad m = 0,\ldots,N_{n}
    \label{lm}.
\end{equation}
Then with some weight functions $w_{m}$, see \cite{trefethen} for a 
discussion and a code to compute them, the 
integral (\ref{gs}) is approximated via
\begin{equation}
     \int_{-1}^{1}g(l)dl\approx \sum_{m=0}^{N_{n}}w_{m}g(l_{m})
    \label{ccalg}.
\end{equation}
Thus for given weights $w_{m}$, $m=0,\ldots,N_{n}$,  this is just a 
scalar product. The Clenshaw-Curtis algorithm is a spectral method, 
for an error analysis see \cite{CC}. 

We are interested in computing the Hilbert transform on the whole 
compactified
real line. For convenience, we use the same discretisation in $x$ as in $y$. Thus infinity becomes a finite point on our numerical grid. 
However, the Hilbert transform is not merely known on the collocation 
points in $x$. For intermediate points we apply a 
numerically stable and efficient interpolation algorithm  in 
the form of \emph{barycentric interpolation}, see \cite{bary} 
for a discussion and references. In this way we obtain the Hilbert 
transform not only at the 
collocation points, but for all $x\in\mathbb{R}\cup\{\infty\}$ we are 
interested in. 

If $x\in I_{n}$, this can lead to an integrand with a 
limit of the type `0/0'. Assuming that the function $f$ is differentiable 
on $I_{n}$, this limit will be calculated via de l'Hospital's rule,
\begin{equation}
    \lim_{x\to y}\frac{f(y)-f(x)}{x-y}=f'(y)
    \label{der}.
\end{equation}

\begin{remark}
	This formula shows also that the terms $\frac{f(y)-f(x)}{x-y}$ 
	appearing in our approach to compute the Hilbert transform are 
	controled in standard way by the derivative $f'(y)$ of the 
	function appearing in the Hilbert transform. Since this 
	derivative is by hypothesis finite, this controls the magnitude 
	of the terms appearing in the quadrature routine. 
\end{remark}
The derivative of $f$ in (\ref{der}) is approximated via Chebyshev differentiation 
matrices, see \cite{trefethen,WR},
\begin{equation}
    \vec{g}^{~\prime}(l) \approx D \vec{g}(l)
    \label{gprime},
\end{equation}
where $\vec{g}$ is the vector with the components 
$g(l_{0}),\ldots,g(l_{N_{n}})$, i.e., $g$ sampled at the Chebyshev 
collocation points. Since $g$ is anyway sampled at these points, it 
is convenient to use a consistent differentiation method. 
For smooth functions and sufficiently small intervals, 
$N_{n}$  can be chosen  small enough so that cancellation errors in the term 
$\frac{f(y)-f(x)}{x-y}$ do not play a role. Note that the limit 
(\ref{gprime}) could also be addressed via a deformation of the 
integration path near $x$ into the complex plane (e.g., a small 
circle). However, since we are interested in applications related to 
dispersive equations as Benjamin-Ono, and since for such equations 
singularities in the complex plane can come close to the real axis, 
see \cite{KR,W2} and references therein, this is not a convenient 
approach in this context.

\subsection{Infinite intervals}\label{infinite}
To treat infinite intervals, we use the local parameter $s=1/y$ for 
$y\sim\infty$. We distinguish two cases, first where the function $f$ 
is analytic in $s$ in a neighborhood of infinity, and second where 
this is not the case. In the first case we consider one interval of the 
form $s\in[\tilde{a},\tilde{b}]$, where $\tilde{a}=1/a_{1}$ and 
$\tilde{b}=1/b_{N_{M-1}}$, in the second case two intervals of 
the form $s\in[\tilde{a},0]$ and $s\in[0,\tilde{b}]$. Thus our 
approach can deal with functions $f$ which do not tend 
to the same finite value for $y\to\pm\infty$ (in the latter case we 
would deal with $N_{M-2}$ finite intervals and two infinite intervals in total). 

In both cases, we get an integral of the form (\ref{I1})
\begin{equation}
   H_{\infty} =  \frac{1}{x}\mathcal{P}\int_{\tilde{a}}^{\tilde{b}}\frac{1}{s}
    f(1/s)\frac{ds}{1/x-s}
    \label{I3}.
\end{equation}
For the following discussion we assume that $f(1/s)/s$ is bounded for 
$s\to0$,  but this is not required in general. Note that we 
discuss in the following section how functions with an essential 
singularity at infinity can be treated. For the remainder of this section we assume that $f$ is 
analytic in $s$ for $s\sim 0$. 

If $1/x\notin [\tilde{a},\tilde{b}]$, the integral (\ref{I3}) can be 
computed as before with the Clenshaw-Curtis algorithm. If 
$1/x\in[\tilde{a},\tilde{b}]$, we proceed as in (\ref{I2}):
\begin{equation}
    H_{\infty} = \frac{1}{x}\mathcal{P}\int_{\tilde{a}}^{\tilde{b}}\frac{1}{s}
    f(1/s)\frac{ds}{1/x-s}=\frac{1}{x}\int_{\tilde{a}}^{\tilde{b}}
    \frac{
    f(1/s)/s-xf(x)}{1/x-s}ds-f(x)\ln\left|
    \frac{1/x-\tilde{b}}{1/x-\tilde{a}}\right|
    \label{I4}.
\end{equation}

The simplest realisation  of our  approach is to compute 
the Hilbert transform on two intervals. For 
convenience we choose here $y\in[-1,1]$ and $1/y\in[-1,1]$.
This leads for (\ref{I1}) to 
\begin{equation}
    \pi\mathcal{H}(f)(x)=\mathcal{P}\int_{-1}^{1}\frac{f(y)}{x-y}dy
    +\mathcal{P}\int_{-1}^{1}\frac{f(1/s)}{s(x/s-1)}ds
    \label{hilbert2}.
\end{equation}
This is equivalent to 
\begin{equation}
    \pi\mathcal{H}(f)(x)=\int_{-1}^{1}\frac{f(y)-f(x)}{x-y}dy
    +\frac{1}{x}\int_{-1}^{1}\frac{f(1/s)/s-xf(x)}{s-1/x}ds
   \label{hilbert3}.
\end{equation}
Thus the logarithmic terms cancel, and we are left with two integrals 
which are defined in a classical sense.

\subsection{Examples}\label{secexample}
We illustrate the above approach with the example of a function which 
is analytic on  $\mathbb{R}\cup\{\infty\}$. Concretely, we consider 
the two functions
\begin{equation}
    f_{1}(y) = \frac{1}{1+y^{2}},\quad f_{2}(y)=\frac{1}{1+y^{4}}
    \label{ex12},
\end{equation}
which are also examples one and two in \cite{weideman}. The Hilbert transforms of 
both functions can be shown to be given in explict form
\begin{equation}
    \mathcal{H}[f_{1}](x) = -\frac{x}{1+x^{2}},\quad 
	\mathcal{H}[f_{2}](x)=
    -\frac{x(1+x^{2})}{\sqrt{2}(1+x^{4})}
    \label{ex12h}.
\end{equation}

We generalize the first example of (\ref{ex12h}) slightly to 
\begin{equation}
    f(y)=\frac{1}{a^{2}+y^{2}},\quad a\in\mathbb{R}^{+}
    \label{example}.
\end{equation}
The Hilbert transform for this function can be calculated with 
(\ref{hilbert3}) to give
\begin{equation}
    \mathcal{H}[f](x) = 
    -\frac{2}{a\pi}(\arctan(1/a)+\arctan a)\frac{x}{a^{2}+x^{2}}   
    \label{example2},
\end{equation}
which gives for $a=1$ the first result in  (\ref{ex12h}).

\begin{figure}[htb!]
 \includegraphics[width=0.49\textwidth]{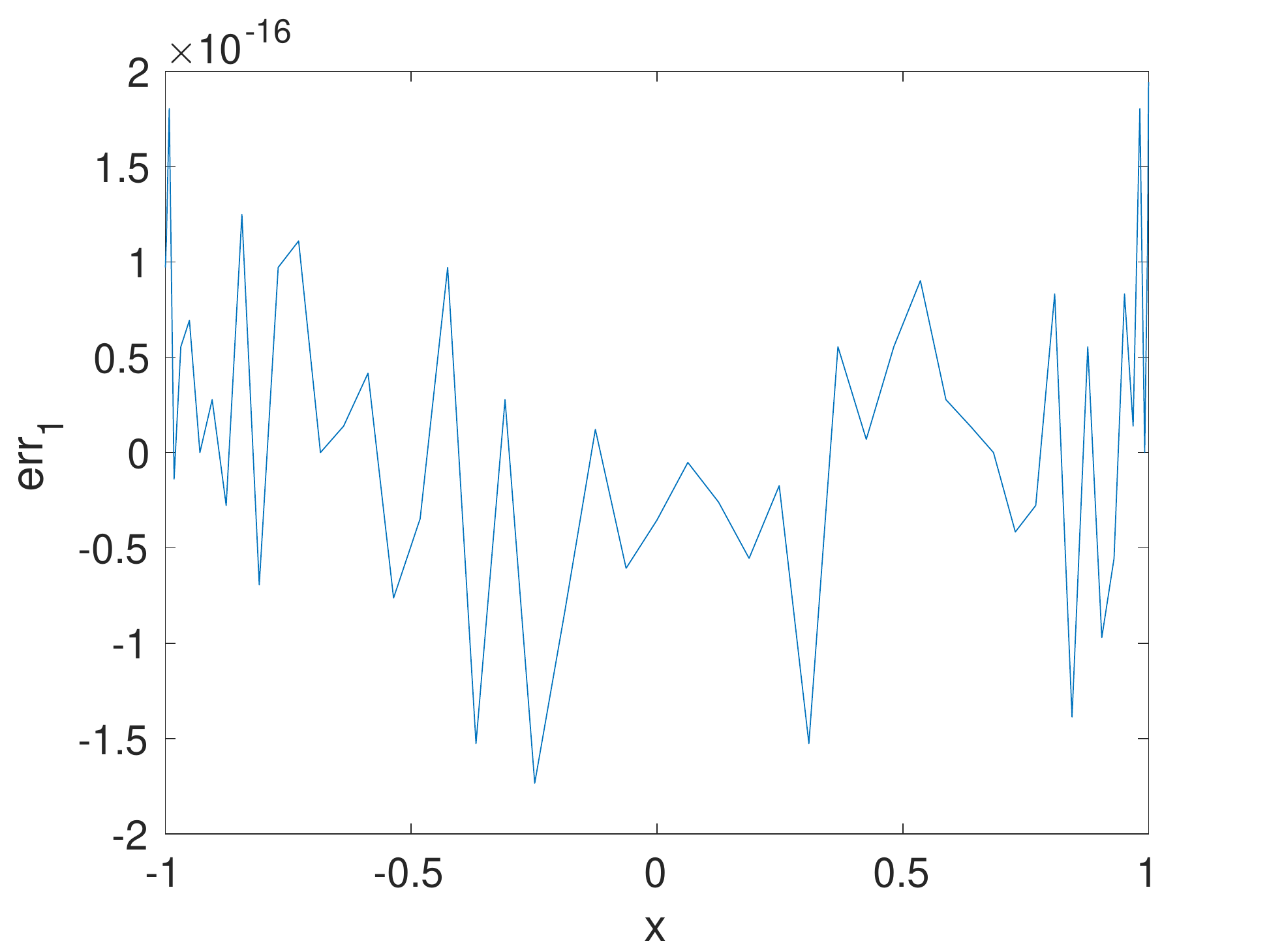}
   \includegraphics[width=0.49\textwidth]{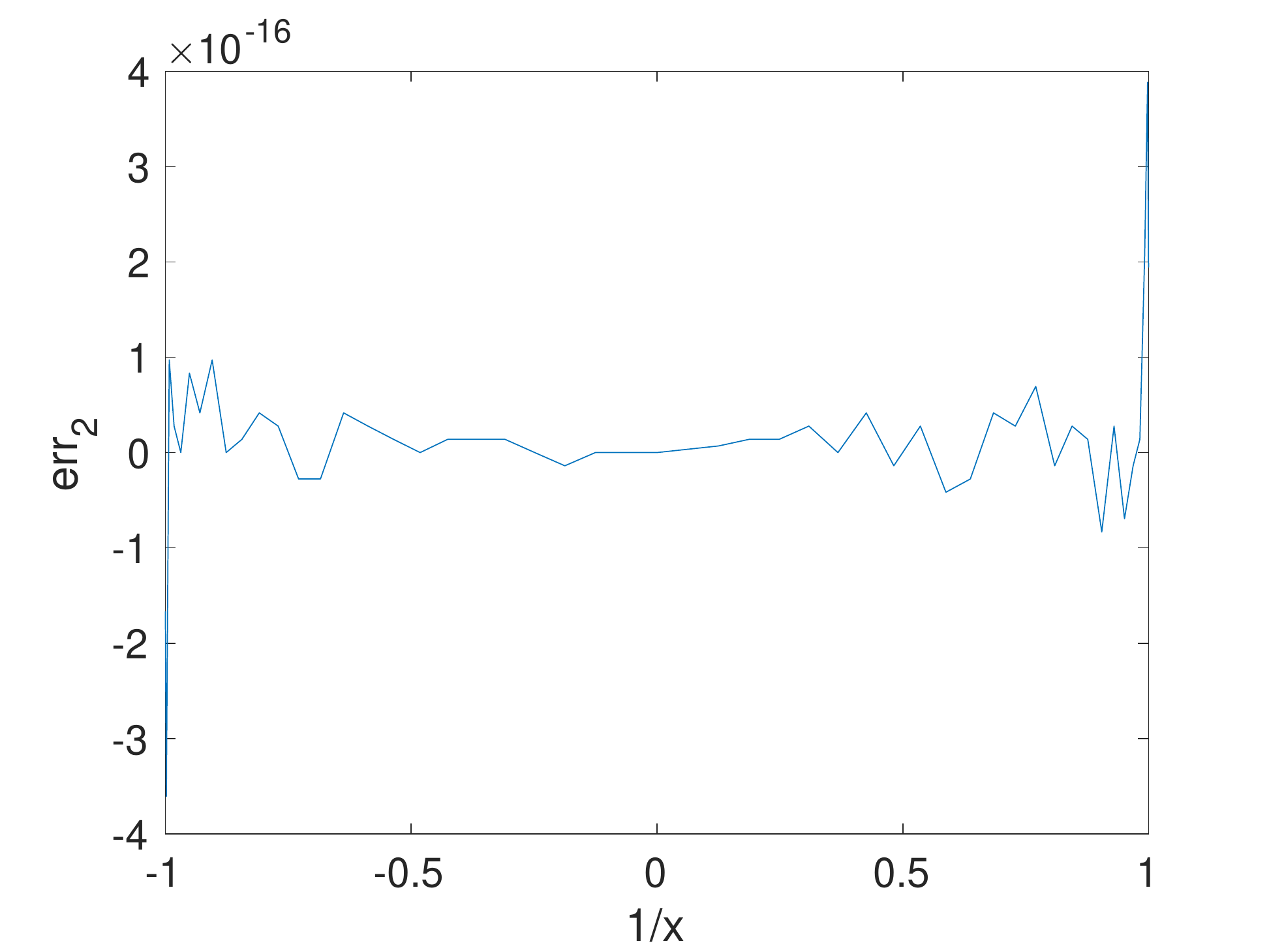}
 \caption{Difference of the first integral in (\ref{hilbert3}) for 
 the function (\ref{example2}) for $a=2$ and $N=50$, on the left for 
 $x\in[-1,1]$, on the right for $1/x\in[-1,1]$.}
 \label{ex1a4err1N50}
\end{figure}

We define as the numerical error $\mbox{err}_{1}$ the difference  
of the first integral in (\ref{hilbert3}) for the function (\ref{example2}) and the 
explicit value  $2/a\arctan(1/a)x/(a^{2}+x^2)$ for $x\in[-1,1]$, and 
$\mbox{err}_{2}$ as the same difference for $1/x\in[-1,1]$. For $N=50$ and $a=4$, 
these errors are shown in Fig.~\ref{ex1a4err1N50}. It can be seen 
that the error is for all values of $x$ of the order of machine 
precision ($10^{-16}$ here). 

\begin{figure}[htb!]
 \includegraphics[width=0.49\textwidth]{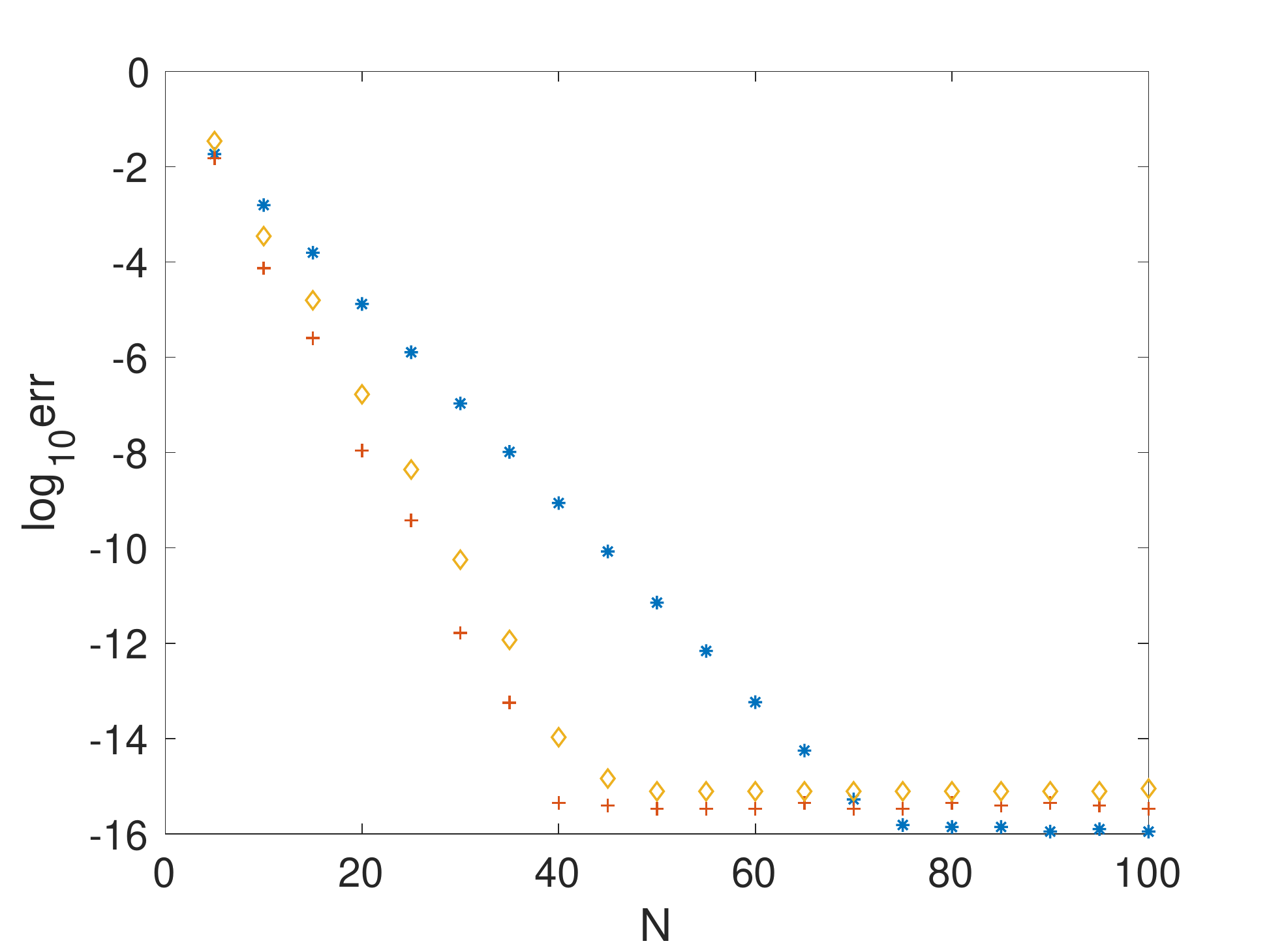}
 \includegraphics[width=0.49\textwidth]{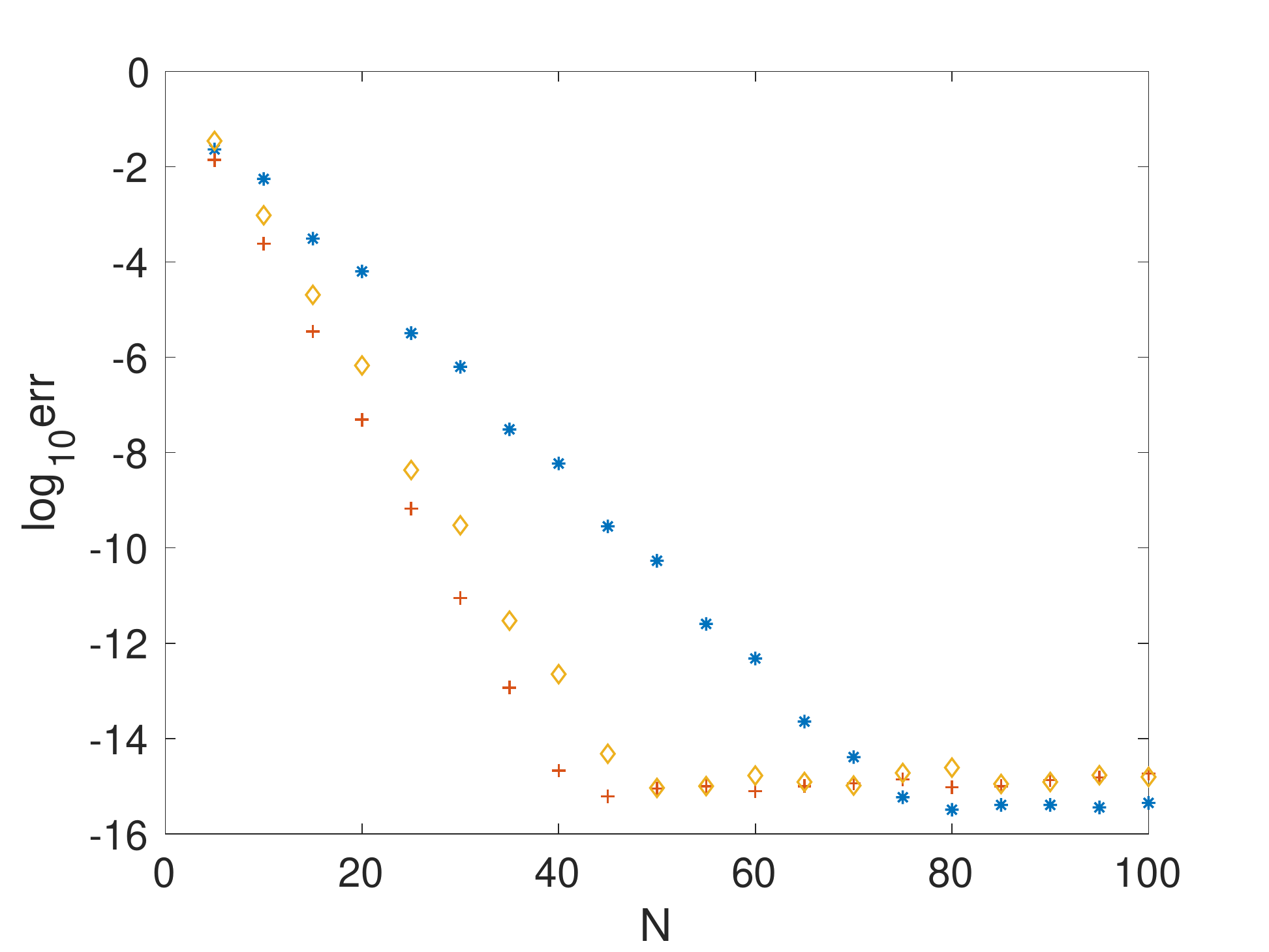}
 \caption{$L^{\infty}$ norm of the difference between the computed Hilbert 
 transform and its exact value in dependence of the number $N_{1,2}$ of 
 collocation points in each of the two intervals; the stars correspond to 
 the example (\ref{example2}) for $a=2$, the plus signs to the same 
 function for $a=1$ (example 1 of \cite{weideman}), and the diamonds 
 to the second case in (\ref{ex12}) (example 2 of \cite{weideman}); on 
 the left with Clenshaw-Curtis, on the right with Gauss-Legendre 
 quadrature.}
 \label{ex12err}
\end{figure}

To study the dependence of the numerical error on 
the number of points in each of the intervals, we define the error 
$\mbox{err}$ as the $L^{\infty}$-norm of the difference between the 
numerically computed Hilbert transform and its exact value in both 
intervals. For 
simplicity we choose the same value of points $N_{1,\infty}$ in both intervals, 
but this is not mandatory. The numerical error can be seen for the 
example (\ref{example2}) for $a=1$ and $a=2$ and for the second 
example in (\ref{ex12}) in Fig.~\ref{ex12err}. The spectral 
convergence of the code can be 
well recognized in a semi-logarithmic plot. The level of the rounding error 
is reached in the cases (\ref{ex12}) for $N_1=40$ points, and for 
(\ref{example2}) with $a=2$ with roughly 70 points. 
\begin{remark}\label{rem2}
    The algorithm discussed above obviously does not depend on the 
    use of Chebyshev collocation points. Instead one could use for 
    instance Gauss-Lobatto points and apply Gauss-Legendre quadrature 
    together with Legendre differentiation matrices for the terms of 
    the form (\ref{der}). The resulting errors are shown on the 
    right of Fig.~\ref{ex12err} for the same examples. As can be seen 
    there is no advantage of the latter algorithm, which is in line 
    with the discussion in \cite{trefethengauss}. We always use 
    Clenshaw-Curtis in the following since we can compute the 
    coefficients of an expansion in terms of Chebyshev polynomials 
    efficiently (see remark \ref{rem3}).
\end{remark}

\begin{figure}[htb!]
 \includegraphics[width=0.49\textwidth]{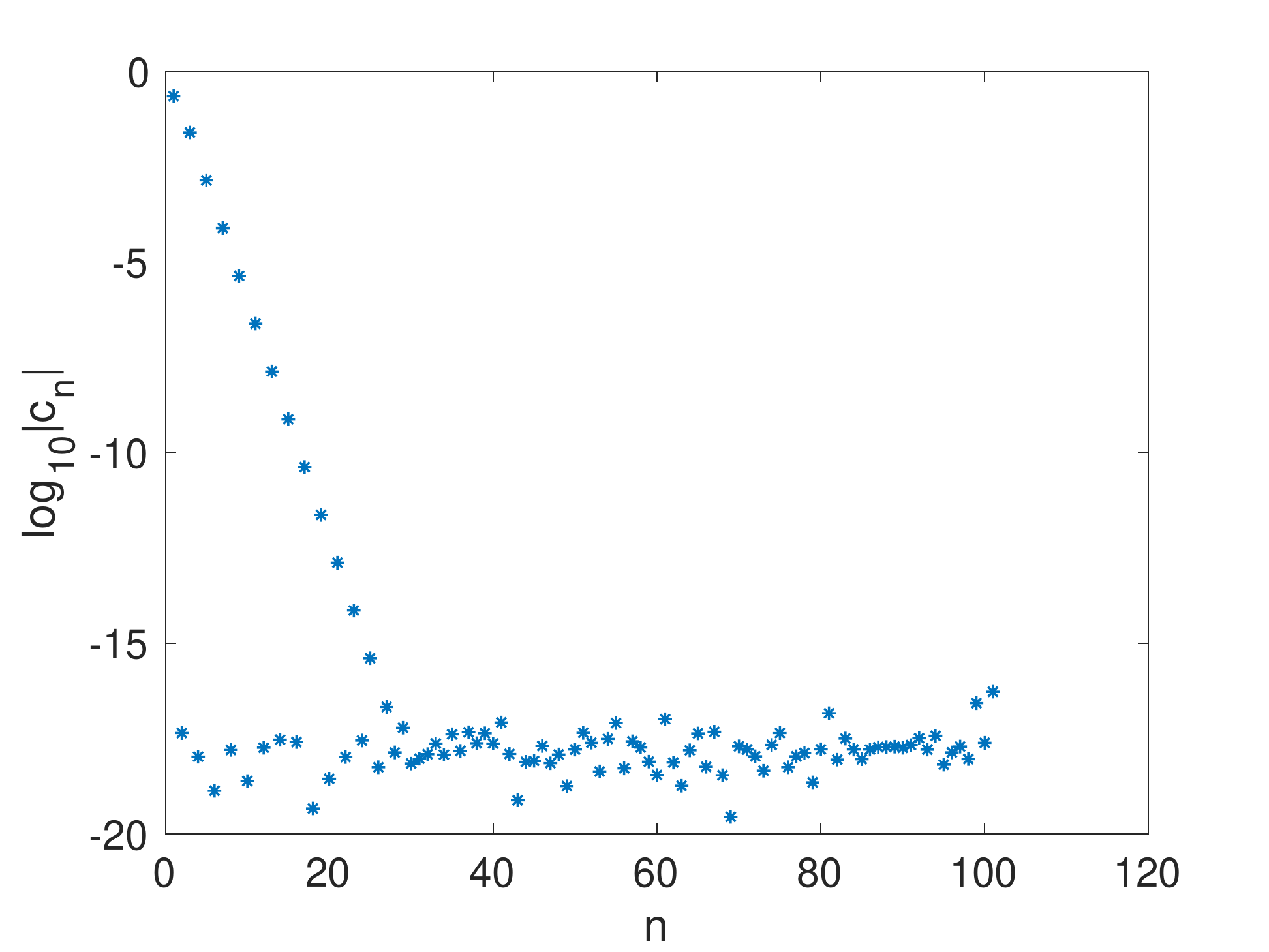}
 \includegraphics[width=0.49\textwidth]{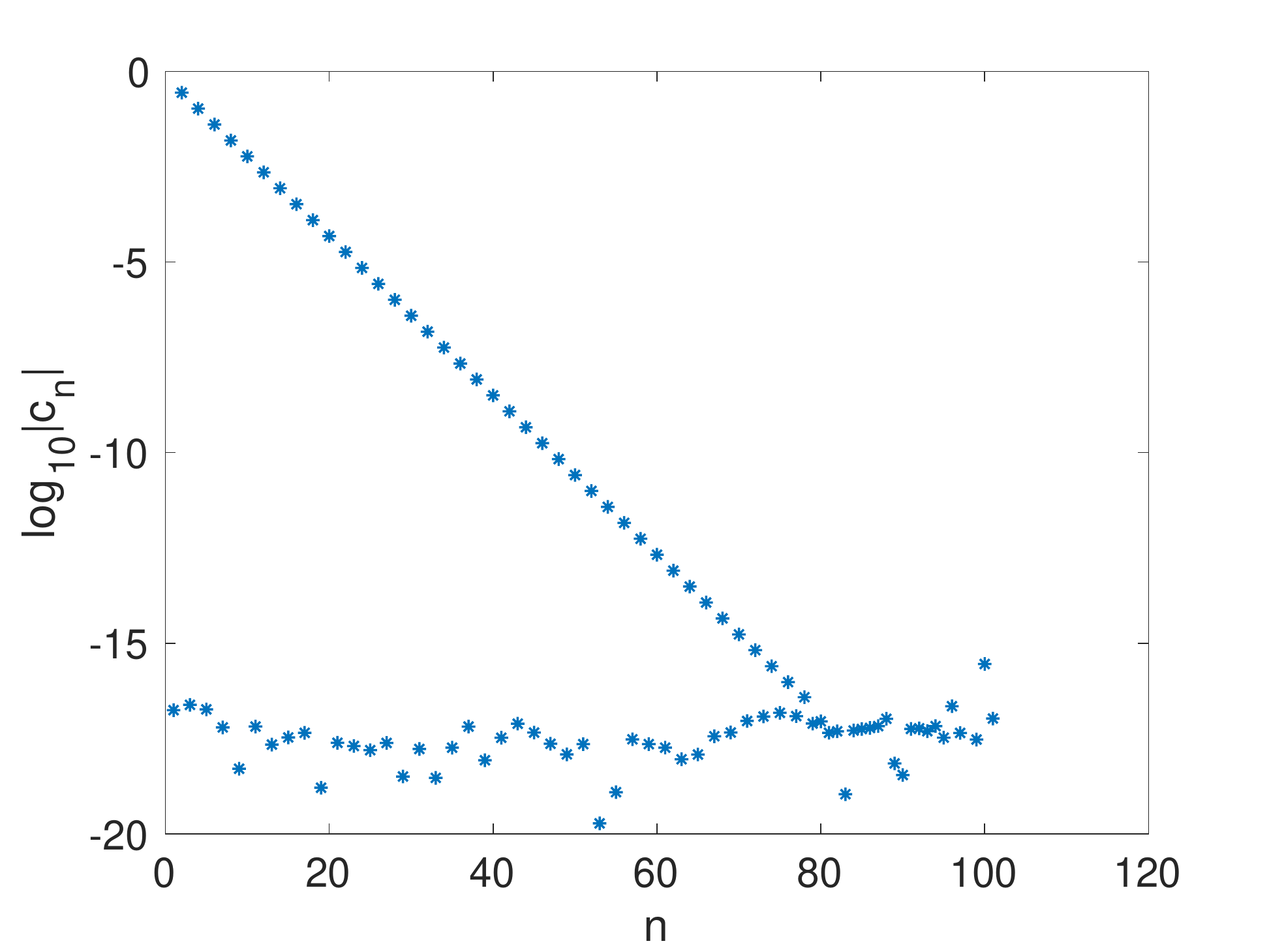}
 \caption{Chebyshev coefficients for the 
 example (\ref{example2}) for $a=2$, on the left for the interval 
$[-1,1]$, on the right for the infinite interval.}
 \label{coeff2}
\end{figure}

\begin{remark}\label{rem3}
    It is evident that the error in the case $a=2$ for function 
    (\ref{example2}) decreases more slowly 
    than in the case $a=1$. Despite this, the error for $y\in[-1,1]$ in 
    Fig.~\ref{ex1a4err1N50} for $N_{1,\infty}=50$ is of the order of machine 
    precision. These two facts indicate that it is not always optimal to 
    choose the same number of points for both intervals. Thus one 
    could either use intervals of the form $x\in[-L,L]$ and 
    $1/x\in[-1/L,1/L]$ with an optimized value of $L$ (this allows to 
	use the same Clenshaw Curtis weights (\ref{ccalg}) in both 
	cases) and the same 
    number of points in both  intervals, or the intervals as 
    before with an optimized value of the number of collocation 
    points for each interval. An indicator for such values can be 
    obtained by considering in each interval the  coefficients  in an 
    approximation of the function $f$ via Chebyshev polynomials \( 
	T_{n}(x):=\cos(n\arccos(x)) \), 
	\begin{equation}
		f(x)
    \approx \sum_{n=0}^{N}c_{n}T_{n}(x)
		\label{chebcoeff}.
	\end{equation}
	These coefficients can be 
    computed efficiently via a \emph{fast cosine transform} which is 
    related to the FFT, see e.g.~\cite{trefethen}. No fast algorithm 
	is known for an expansion in terms of Legendre polynomials. 
\end{remark}
 For the example 
    (\ref{example2}), these Chebyshev coefficients are shown in 
    Fig.~\ref{coeff2}. It can be seen that $N_1 = 30$ points are sufficient 
    on the finite interval to reach machine precision, whereas more 
    than $N_{\infty} = 80$ are needed on the infinite interval. 

\subsection{Weideman's approach to the computation of the Hilbert 
transform on the real line}
As mentioned in the introduction, Weideman's approach \cite{weideman} 
is based on the mapping of the real line to the circle,
\begin{equation}
    y = \tan \frac{\theta}{2},\quad \theta\in[-\pi,\pi].
    \label{map}
\end{equation}
His approach uses an expansion of the considered functions 
in terms of rational functions
\begin{equation}
    \phi_{n} = \frac{(1+iy)^{n}}{(1-iy)^{n+1}},\quad n\in\mathbb{Z}
    \label{phin},
\end{equation}
instead of trigonometric ones: with (\ref{map}) one 
obviously has that 
\begin{equation}
    f(y)=\sum_{n\in\mathbb{Z}}^{}a_{n}\phi_{n}(y)\Rightarrow 
    f(y)(1-iy)=\sum_{n\in\mathbb{Z}}^{}a_{n}e^{in\theta}
    \label{fphi}.
\end{equation}
The Hilbert transform acts on the $\phi_{n}$ as $\mathcal{H}\phi_{n}=
i\mbox{sgn}(n)\phi_{n}$, $n\in \mathbb{Z}$, see \cite{weideman}.  
On the latter an FFT approach is implemented,
\begin{equation}
	f(y)(1-iy)\approx \sum_{-N_{\mathcal{F}}/2}^{N_{\mathcal{F}}/2-1}a_{n}e^{in\theta}
	\label{fapprox},
\end{equation}
where \( N_{\mathcal{F}} \) is an even natural number, and where the coefficients 
\( a_{n} \), \( n=-N_{\mathcal{F}}/2,\ldots,N_{\mathcal{F}}/2-1 \) are computed with an FFT. The 
Hilbert transform is thus approximated as 
\begin{equation}
	\mathcal{H}(f)\approx \frac{1}{1-ix}\sum_{n=-N_{\mathcal{F}}/2}^{N_{\mathcal{F}}/2-1}i\mbox{sgn}(n)
	a_{n}e^{in\theta},
	\label{happrox}
\end{equation}
with the definition \( \mbox{sgn}(0)=1 \). To approximate the Hilbert 
transform in this way, two FFTs are necessary. Note that we use 
here the \( N_{\mathcal{F}} \) as e.g.\ in \cite{trefethen} for the FFT, which is 
twice the value used in \cite{weideman}. Since Weideman 
\cite{weideman} compared his method to several numerical approaches, 
we will only relate our results to his in this paper. 

\begin{figure}[htb!]
 \includegraphics[width=0.49\textwidth]{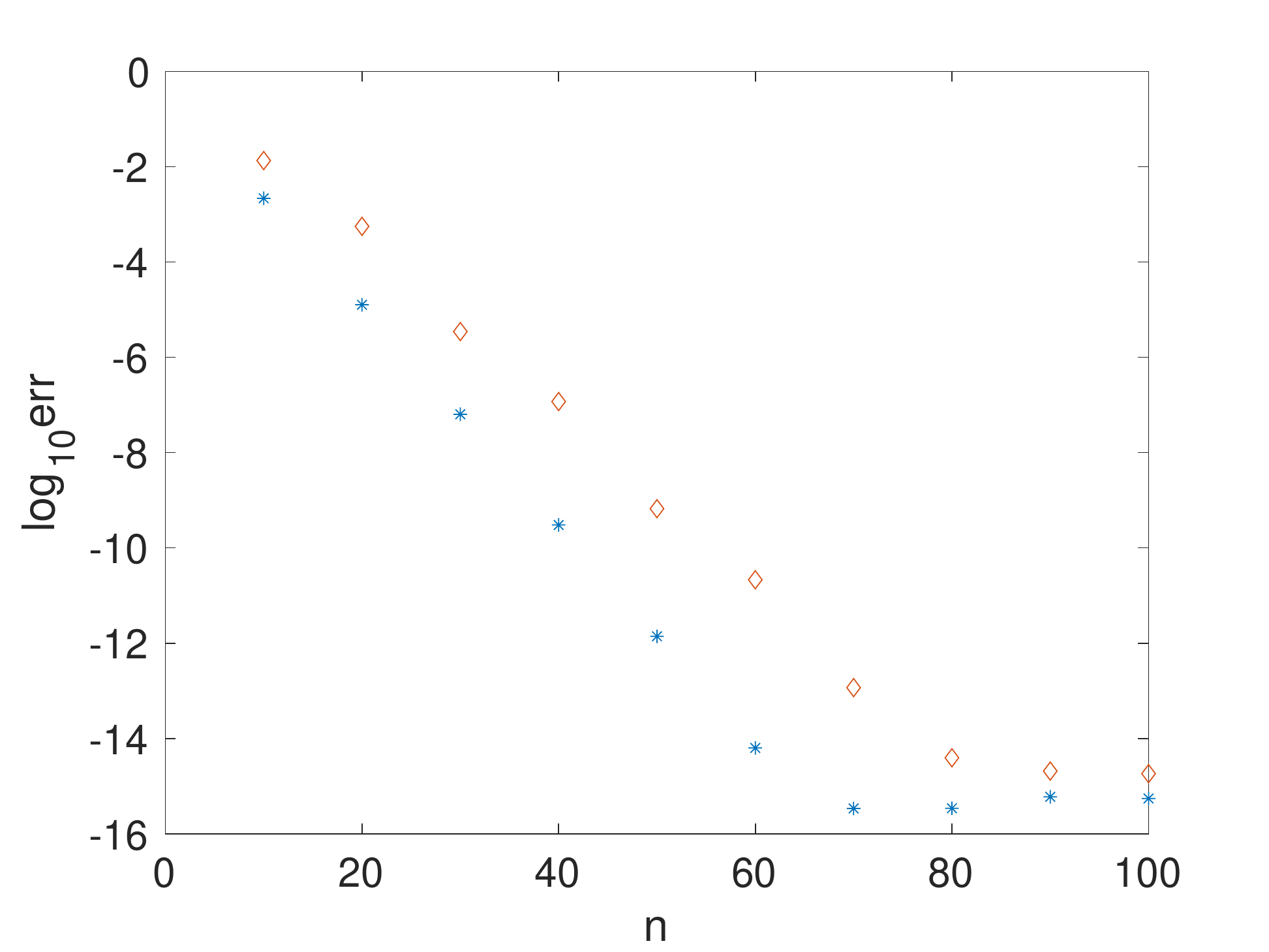}
 \includegraphics[width=0.49\textwidth]{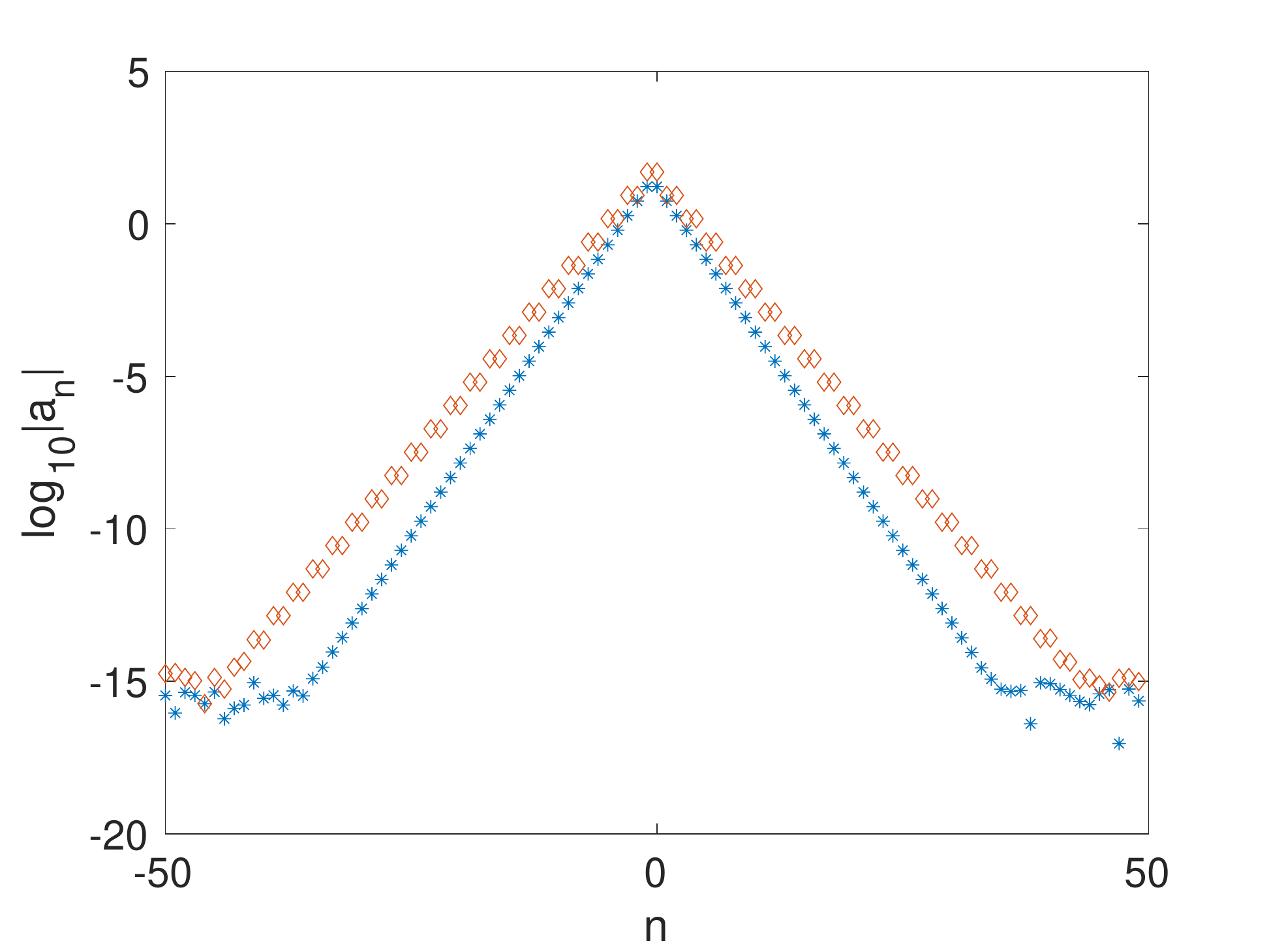}
 \caption{$L^{\infty}$-norm of the difference between the computed Hilbert 
 transform (with the method \cite{weideman}) and its exact value in 
 dependence of the number $n$ of  collocation points on  the left; the stars correspond to  the example (\ref{example2}) for $a=2$,  and the diamonds 
 to the second case in \ref{ex12}; on the right the coefficients \( 
 a_{n} \) for both examples.}
 \label{AWerror12}
\end{figure}

The first example in (\ref{ex12h}) is trivial in this approach since 
\( f_{1}(y)=(\phi_{0}(y)+\phi_{-1}(y))/2 \) which also gives the 
formula for the Hilbert transform. The numerical errors in dependence 
of the number \( N_{\mathcal{F}} \) in (\ref{fapprox}) for the second example in 
(\ref{ex12h}) and (\ref{example2}) are shown in 
Fig.~\ref{AWerror12} on the left. Spectral convergence is evident. The approach reaches machine precision with roughly the 
same number of points as used above in each of the two intervals, 
i.e., half the number of collocation points in total (an optimal choice would 
be \( N_{1}+N_{\infty}\sim 110 \) compared to 80 points in Weideman's 
approach). The values of \( N_{\mathcal{F}} \) 
needed to reach machine precision can be as usual estimated via the 
coefficients \( a_{n} \), \( n=-N_{\mathcal{F}}/2,\ldots,N_{\mathcal{F}}/2-1 \) which are shown 
on the right of Fig.~\ref{AWerror12}. 

It is thus not a  surprise that the approach \cite{weideman} is somewhat 
more efficient for functions analytic on the whole real line. Still, the 
order of magnitude of the total number of collocation points  to achieve machine precision is the same as for the Weideman approach and the
multi-domain approach in the present paper.

\section{Piece-wise analytic functions}
Multi-domain spectral methods are especially efficient for functions 
which are not globally smooth, and this will be illustrated in the 
present section. We will address the example of two 
intervals as in the previous section, but with functions which are 
only continuous or not even that. Since the logarithms in formulae 
(\ref{I2}) and (\ref{I4}) are taken care of analytically, only the 
integrals there have to be computed numerically. These integrals have smooth integrands and can be efficiently computed. Note 
that the logarithms will lead to unbounded terms for functions which 
are not continuous on \( \mathbb{R} \).

\begin{remark}
	For values of \( x \) in the second interval, the integral in 
	(\ref{I2}) over the first interval can be computed as in 
	(\ref{I1}). But for \( x \) close to the boundary of the 
	interval, this leads to an almost singular integrand which is 
	difficult to approximate with polynomials. This is why we insist 
	here on piecewise analytic functions which allow for an analytic 
	continuation of the function in the interval \( I_{1} \) to a 
	slightly larger interval. The second line of (\ref{I2}) is used 
	with this analytic continuation for \( x \) close to the 
	boundaries. In this way the integrand is always controlled by the 
	derivative of $f$. 
\end{remark}

\begin{figure}[htb!]
 \includegraphics[width=0.49\textwidth]{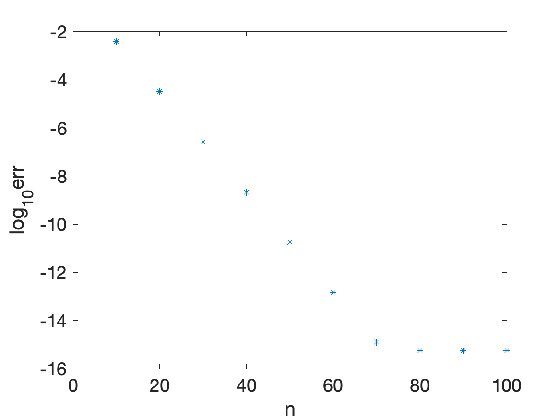}
 \includegraphics[width=0.49\textwidth]{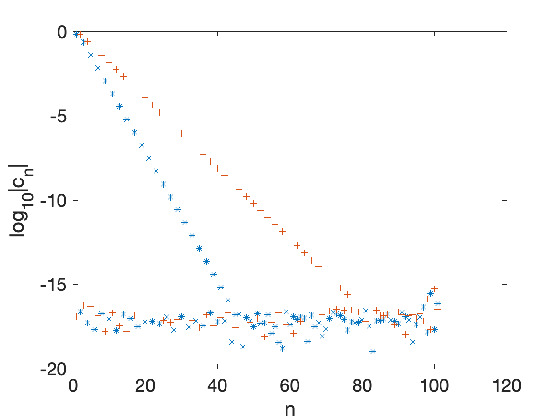}
 \caption{$L^{\infty}$ norm of the difference between the computed Hilbert 
 transform  and its exact value (\ref{Hf2int}) in dependence of the number $N$ of 
 collocation points in each domain on 
 the left;  on the right the Chebyshev coefficients \( 
 c_{n} \) for both examples, the stars corresponding to the finite 
 interval, the plus signs to the infinite one.}
 \label{fconterr}
\end{figure}

Concretely we will address the 
example
\begin{equation}
	f = 
	\begin{cases}
		\frac{1}{a_{1}^{2}+y^{2}}, & |y|\leq 1 \\
		\frac{\alpha}{a_{2}^{2}+y^{2}}, & |y|>1
	\end{cases}
	\label{f2int},
\end{equation}
where \( a_{1} \), \( a_{2} \), \( \alpha \) are constants. Each 
function in the respective interval 
has an obvious analytic continuation to the whole real line. 
The 
Hilbert transform of (\ref{f2int}) is with relations (\ref{I2}) and 
(\ref{I4})

\begin{equation}
\begin{array}{r}
	\mathcal{H}[f](x) = 
	\frac{2}{a_{1}}\arctan(1/a_{1})\frac{x}{a_{1}^2+x^2}+\frac{2\alpha}{a_{2}}\arctan(a_{2})\frac{x}{a_{2}^2+x^2} \hspace{2cm} \\
	-\left(\frac{1}{a_{1}^2+x^2}-\frac{\alpha}{a_{2}^2+x^2}\right)\ln\left|\frac{1-x}{1+x}\right|  . 
	\label{Hf2int}
\end{array}
\end{equation}
~\\
We consider first the case of a continuous potential, \( 
\alpha=(a_{2}^{2}+1)/(a_{1}^{2}+1) \), where the Hilbert 
transform is bounded, and we choose \( a_{1}=1 \) and \( a_{2}=2 \). 
We use again the same number \( N_{1,\infty} \) of collocation points in both 
intervals. The numerical error (as before, the \( L^{\infty} \) norm 
of the difference between numerical and exact solution) in dependence 
of \( N_{1,\infty} \) can be seen on the left of Fig.~\ref{fconterr}. As expected 
the error decreases exponentially and reaches machine precision at 
essentially the same values as in the previous section. This means 
that as theoretically predicted, only the regularity on the 
respective intervals is important. As in the previous section, it is 
not optimal to choose the same resolution in both intervals. This is 
indicated by the decrease of the Chebyshev coefficients on the right 
of Fig.~\ref{fconterr} for \( N_{1,\infty}=100 \). They decrease in both cases 
exponentially, but more rapidly in the finite domain. Thus as in the 
examples of the previous section, \( N_{1}+N_{\infty}\sim120 \) would 
allow to achieve machine precision with two domains. 

\begin{figure}[htb!]
 \includegraphics[width=0.32\textwidth]{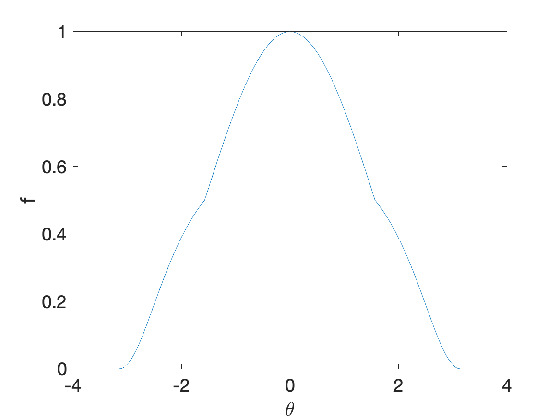}
 \includegraphics[width=0.32\textwidth]{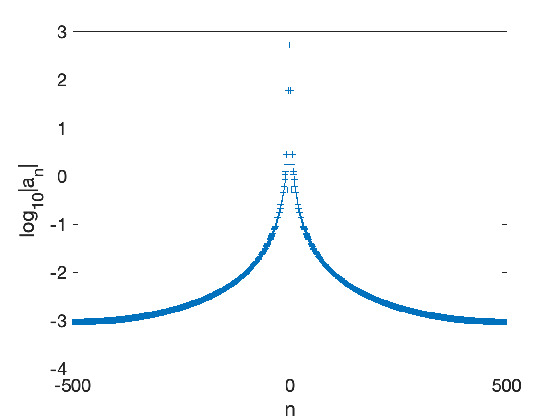}
 \includegraphics[width=0.32\textwidth]{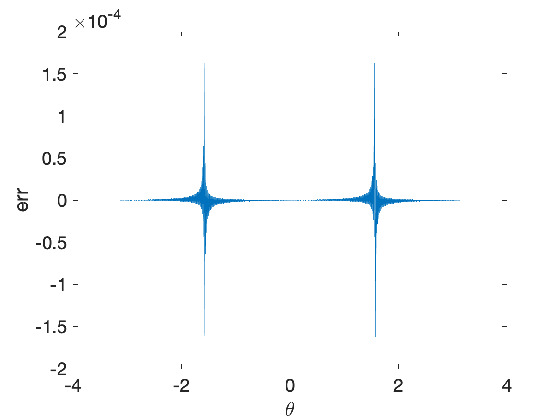}
 \caption{The function (\ref{f2int}) in dependence of \( \theta \) on 
 the left, its coefficents \( a_{n} \) in dependence of \( n \) for 
 \( N_{\mathcal{F}}=1000 \) in the middle, and the $L^{\infty}$ norm of the difference between the computed Hilbert 
 transform (with the method \cite{weideman}) and its exact value 
 (\ref{Hf2int}) for $N_{\mathcal{F}}=1000$  on 
 the right.}
 \label{AWconterr}
\end{figure}

The situation is very different for the global approach 
\cite{weideman} for which only the regularity on the whole 
compactified real line counts. The function (\ref{f2int}) in 
dependence of the coordinate \( \theta \) on the circle to which the 
real line is mapped can be seen on the left of Fig.~\ref{AWconterr}. 
The corresponding coefficients \( a_{n} \) for this function can be 
seen in the same figure in the middle. As expected for a piecewise 
continuous function, they only decrease algebraically, for \( N_{\mathcal{F}}=1000 
\) only to the order of \( 10^{-3} \). The difference of the numerical 
and the exact solution for \( N_{\mathcal{F}}=1000 \) can be seen on the right of 
Fig.~\ref{AWconterr}. It is of the order of \( 10^{-4} \) where the 
main error comes as expected from the domain boundaries where the 
function is not differentiable. 

\begin{figure}[htb!]
 \includegraphics[width=0.49\textwidth]{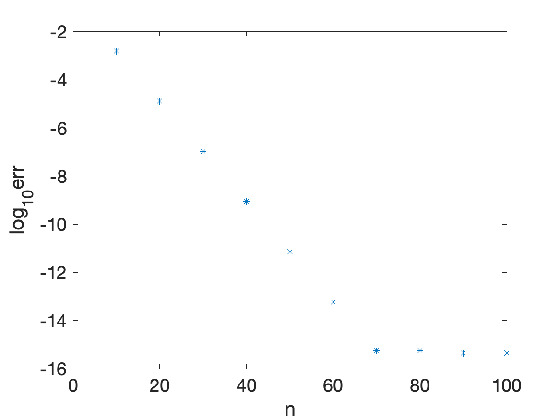}
 \includegraphics[width=0.49\textwidth]{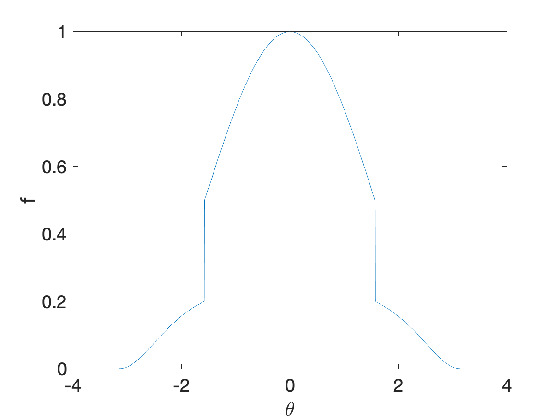}
 \caption{$L^{\infty}$ norm of the difference between the computed Hilbert 
 transform  and its exact value (\ref{Hf2int}) in dependence of the number $N_{1,\infty}$ of 
 collocation points in each domain on 
 the left for the function shown on the right;  on the right the 
 function (\ref{f2int}) in dependence of the coordinate \( \theta \) 
 for \( a_{1}=1 \), \( a_{2}=2 \) and \( \alpha=1 \).}
 \label{fdisconterr}
\end{figure}

For discontinuous potentials being analytic on the respective 
intervals, not much changes for the multi-domain approach. If we 
consider the same example as in Fig.~\ref{fconterr}, just with \( 
\alpha=1 \), the error on the left of Fig.~\ref{fdisconterr} shows 
virtually the same behavior as in Fig.~\ref{fconterr}. This is due 
to the fact the Chebyshev coefficients are the same up to 
multiplication by the factor \( \alpha \). As a function on the whole 
real line, \( f \) is now obviously discontinuous, see the figure on 
the right of Fig.~\ref{fdisconterr}.

\begin{figure}[htb!]
 \includegraphics[width=0.49\textwidth]{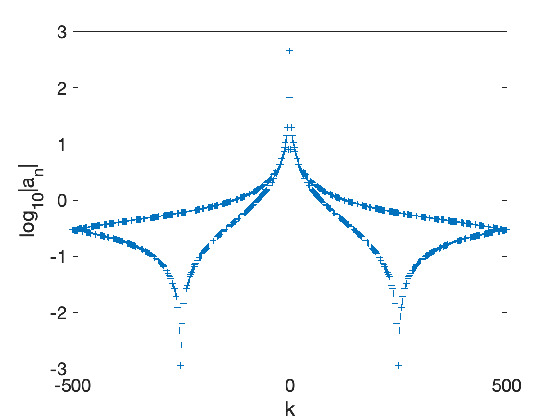}
 \includegraphics[width=0.49\textwidth]{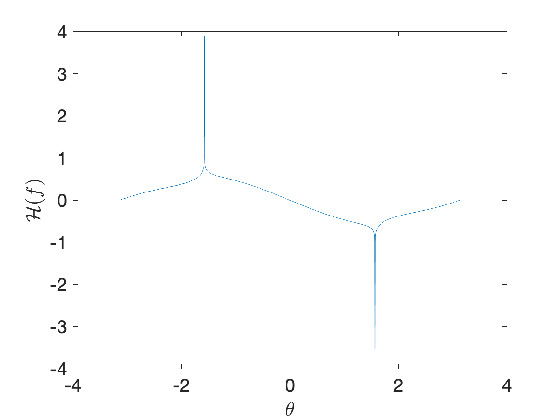}
 \caption{The coefficients \( a_{n} \) for the function on the right 
 of Fig.~\ref{fdisconterr} for \( N_{\mathcal{F}}=1000 \) on the left, and the Hilbert 
 transform (\ref{Hf2int}) on the right.}
 \label{fdiscontcoeff}
\end{figure}

The discontinuity of \( f \) implies that the global approach 
\cite{weideman} leads to a Gibbs phenomenon at the discontinuities, 
and the coefficients \( a_{n} \) consequently decrease only very 
slowly, see the left of Fig.~\ref{fdiscontcoeff}. The situation for 
the Hilbert transform is worse since the latter has a logarithmic 
divergence at the discontinuities as shown on the right of 
Fig.~\ref{fdiscontcoeff} (the values where the logarithm becomes 
infinite are obviously not shown). As is well known, logarithms are not 
efficiently approximated by Fourier series. 

\section{Essential singularities at infinity}
The focus of this paper is on the Hilbert transform of functions 
which are piecewise analytic on the compactified real axis. As has 
been shown in the previous sections for various examples, spectral convergence is achieved 
in such cases. For completeness we add here the remaining examples of 
\cite{weideman} which have essential singularities at infinity. The 
polynomial methods applied in the present paper are not ideal in such 
a case, but as we will show in this section, can still be used 
successfully. We first discuss the case of rapidly decreasing 
functions where the integration is just performed on  finite 
intervals. We also add the case of oscillatory singularities at 
infinity for which deformation techniques are applied. 

\subsection{Rapidly decreasing functions}
For rapidly decreasing functions, the main change with respect to the 
previous sections is that no integration on an infinite interval is 
needed. In the simplest case one just works on \( [-L,L] \) where \( 
L>0 \) is chosen such that \( f \) vanishes with numerical precision 
for \( |x|>L \). 

The first example in this context is example 5 of \cite{weideman}, 
the Gauss function with the Hilbert transform
\begin{equation}
	\mathcal{H}\left(e^{-y^{2}}\right)=-\frac{2}{\sqrt{\pi}}D(x),\quad
	D(x) = e^{-x^{2}}\int_{0}^{x}e^{t^{2}}dt;
	\label{dawson}
\end{equation}
here \( D(x) \) is Dawson's integral which we compute with the 
corresponding function in Octave (no tolerance is given there, but the 
 results below indicate it is computed with machine precision). 

\begin{figure}[htb!]
 \includegraphics[width=0.32\textwidth]{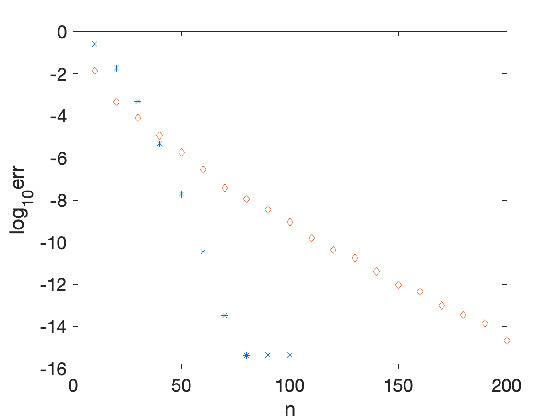}
 \includegraphics[width=0.32\textwidth]{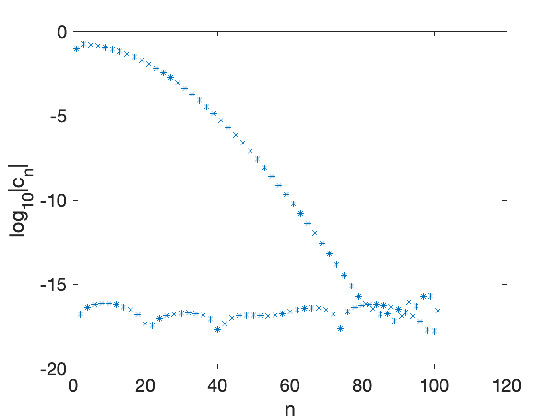}
 \includegraphics[width=0.32\textwidth]{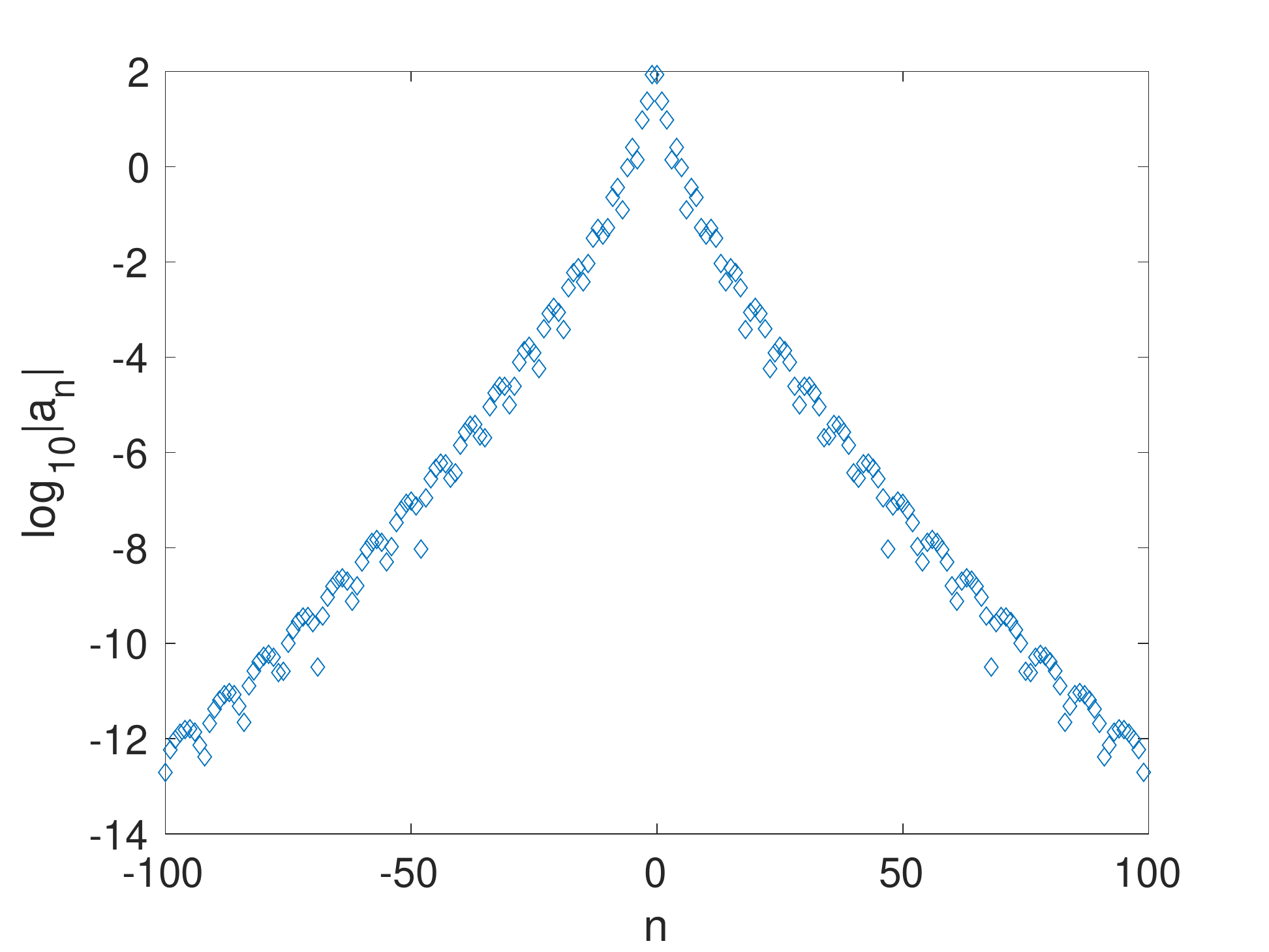}
 \caption{The numerical error for the computation of the Hilbert 
 transform for the Gaussian on the left (`stars' for the multi-domain 
 method, `diamonds' for the global approach), the Chebyshev coefficients
  for the Gaussian in the middle, and its coefficients \( a_{n} \) in dependence of \( n \) for 
 \( N=1000 \)  on 
 the right.}
 \label{dawsonerr}
\end{figure}

To compute the Hilbert transform of the Gauss function, we choose \( 
L = 6 \) (as usual in a way that the spectral coefficients decrease 
exponentially). The numerical error in the computation of the Hilbert 
transform can be seen on the left of Fig.~\ref{dawsonerr}. The error 
for the multi-domain approach (stars) 
decreases exponentially as expected and reaches machine precision 
with roughly 80 collocation points. In the same figure we show 
(diamonds) the corresponding error for the global approach 
\cite{weideman}. Here around $N_{\mathcal{F}} = 200$ collocation points are needed to 
reach the same precision. This is obviously due to the essential 
singularity at infinity which is simply omitted (we work on a finite 
interval) in the multi-domain approach, but which is important in the 
global approach on the compactified real axis. This can be also seen 
from the spectral coefficients, in the middle of Fig.~\ref{dawsonerr} 
for the Chebyshev coefficients of the multi-domain approach, and on 
the right of the same figure the coefficients of the expansion in 
terms of rational functions. 

\begin{figure}[htb!]
 \includegraphics[width=0.32\textwidth]{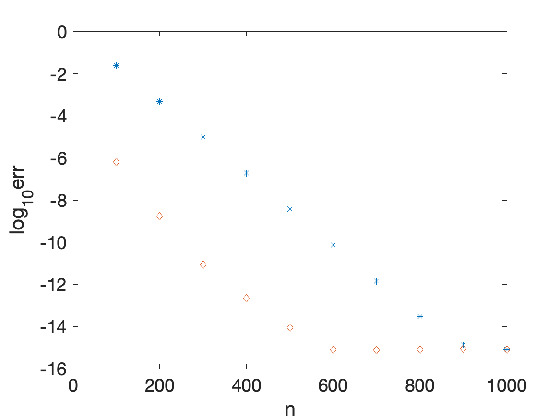}
 \includegraphics[width=0.32\textwidth]{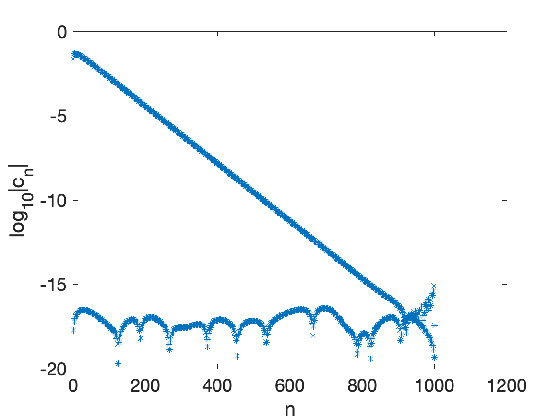}
 \includegraphics[width=0.32\textwidth]{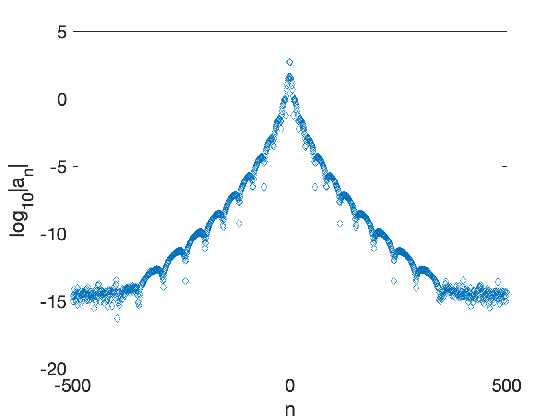}
 \caption{The numerical error for the computation of the Hilbert 
 transform for the hyperbolic secans on the left (`stars' for the multi-domain 
 method, `diamonds' for the global approach), 
 the Chebyshev coefficients
  for the function in the middle, and its coefficients \( a_{n} \) in dependence of \( n \) for 
 \( N=1000 \)  on 
 the right.}
 \label{secherr}
\end{figure}

The situation changes somewhat if the function is just exponentially  
decreasing towards infinity, i.e., if the decrease is 
slower than for the Gaussian. Example 6 of \cite{weideman} is the 
function \( \mbox{sech } x \) for 
which the Hilbert transform is given by 
\begin{equation}
	\mathcal{H}[\mbox{sech}](x)= \tanh(x) + \frac{i}{\pi}\left[
	\psi\left(\frac{1}{4}+\frac{ix}{2\pi}\right)- 
	\psi\left(\frac{1}{4}-\frac{ix}{2\pi}\right) \right]
	\label{sech},
\end{equation}
where the digamma function \( \psi \) is given by the logarithmic 
derivative of the gamma function, \( \psi(z)=\partial_{z}\ln 
\Gamma(z) \). For the multi-domain approach we again use only one 
interval which has to be much larger (\( [-40,40] \)) here because of 
the slower decay of the hyperbolic secans for \( |y|\to\infty \) than 
the Gaussian. This also implies that with both the global and the 
multi-domain approach much higher resolutions are needed than in 
Fig.~\ref{dawsonerr}. The numerical error is shown on the right of 
Fig.~\ref{secherr}. The global approach \cite{weideman} reaches 
machine precision for roughly $N_{\mathcal{F}} = 600$ collocation points, the 
multi-domain approach for roughly $N = 900$ points. This is in accordance 
with the spectral coefficients shown in the same figure in the middle 
and on the right respectively. This implies that in contrast to the 
case of the Gaussian, the global approach for this example is  more efficient than the 
multi-domain approach. 

Example 7 of \cite{weideman} is the function \( f(x) = \exp(-|x|) \), 
a rapidly decreasing function which is smooth on \( \mathbb{R}^{\pm} 
\), but not on \( \mathbb{R} \) and thus an interesting test for a 
multi-domain approach. The Hilbert transform of this function reads
\begin{equation}
	\mathcal{H}[e^{-|y|}] = 
	-\frac{1}{\pi}\mbox{sgn}(x)\left(e^{|x|}E_{1}(|x|)+e^{-|x|}Ei(|x|)\right)
	\label{Heabs},
\end{equation}
where 
\[ E_{1}(x) = \int_{x}^{\infty}e^{-t}\frac{dt}{t},\quad Ei(x) = 
-\mathcal{P}\int_{-x}^{\infty}e^{-x}\frac{dt}{t},
\]
i.e., \( E_{1}(-x)=-Ei(x)-i\pi \).
For the multi-domain approach we use here 3 intervals, \( [-L,0] \), \( 
[0,L] \) and the infinite interval \( |x|>L \). If we choose, as for 
the example in Fig.~\ref{secherr} \( L=40 \), the integral over the 
third interval vanishes with numerical accuracy. The integral on the 
interval \( [-L,0] \) for values of \( x\in[0,L] \) can be 
regularized for \( x\sim0 \) as in (\ref{I2}). This is problematic for some 
values of \( x\in[0,L] \) since \( \exp(x) \) is exponentially growing 
in this interval. Therefore we use the regularization only for 
values of \( x<L_{0} \) with \( L_{0}\sim 1 \) (we take \( L_{0}=1 \) 
here), and similarly for the integral over the interval \( [0,L] \) 
for \( x\in[-L,0] \). 

\begin{figure}[htb!]
 \includegraphics[width=0.32\textwidth]{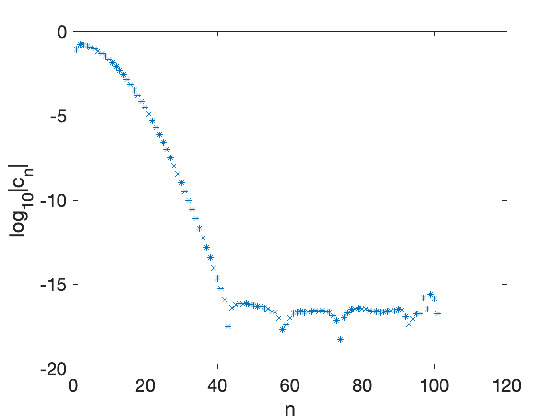}
 \includegraphics[width=0.32\textwidth]{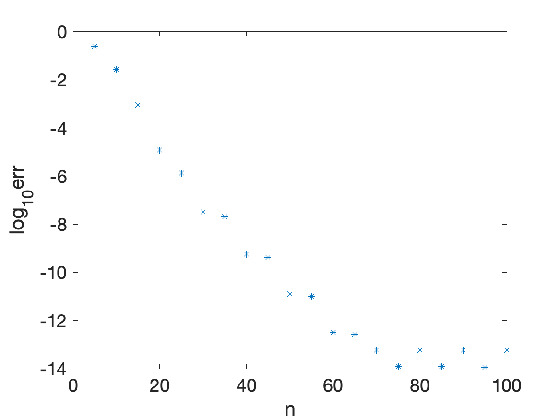}
 \includegraphics[width=0.32\textwidth]{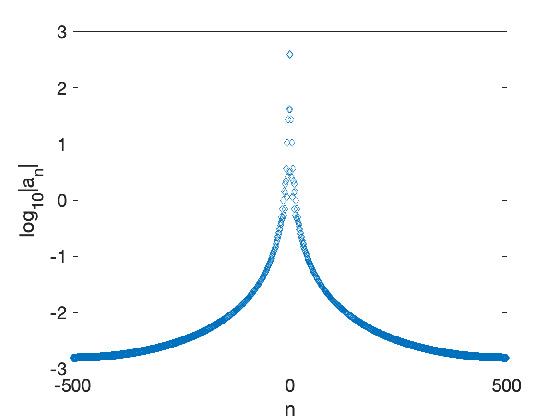}
 \caption{The numerical error for the computation of the Hilbert 
 transform for the function \( f(x)=e^{-|x|} \) on the left, the 
 Chebyshev coefficients
  for the function in the interval \( [-40,0] \) in the middle, 
  and its coefficents \( a_{n} \) in dependence of \( n \) for 
 \( N_{\mathcal{F}}=1000 \)  on 
 the right.}
 \label{abserr}
\end{figure}

The spectral coefficients for the multi-domain approach can be seen  on 
the left of Fig.~\ref{abserr} (we show only the coefficients for \( 
x\in[-L,0] \) for symmetry reasons), and for the global approach 
\cite{weideman} in the middle of the same figure.  
Machine precision is reached in the former case with just 40 
collocation points, whereas in the latter the coefficients for \( 
N=10^{3} \) decrease to the order of \( 10^{-3} \). The numerical 
error for the multi-domain approach can be seen on the right of 
Fig.~\ref{abserr}\footnote{The computation of the exact solution 
(\ref{Heabs}) is not trivial since the function \( Ei(x) \) computed 
in Octave via the function \emph{expint} has to be multiplied with an 
exponential. Thus we write (\ref{Heabs}) in the form \( \mathcal{H}(e^{-|x|}) = 
	-\frac{1}{\pi}\mbox{sgn}(x)\left(Y(x)-Y(-x)\right)\), where \( 
	Y=e^{x}Ei(x) \). The function \( Y \) satisfies the differential 
	equation \( Y'+Y=1/x \) which is solved as in \cite{CFKSV} on 3 
	intervals \( [-\infty,-1] \), \( [-1,1] \), \( [1,\infty] \) with 
	the asymptotic condition \( \lim_{x\to-\infty}Y=0 \) and the 
	substitution \( Y=\tilde{Y}+\ln x e^{-x} \) in the finite interval. As shown in 
	\cite{CFKSV}, machine precision can be reached with such a hybrid
	approach. 
 }. It can be recognized that machine precision is reached with \( 
 N\sim 70\) (the error for the global approach \cite{weideman} for \( 
 N_{\mathcal{F}}=1000 \) is of the order of \( 10^{-4} \)). 

\subsection{Oscillatory singularities at infinity}
\begin{figure}[htb!]
 \includegraphics[width=0.32\textwidth]{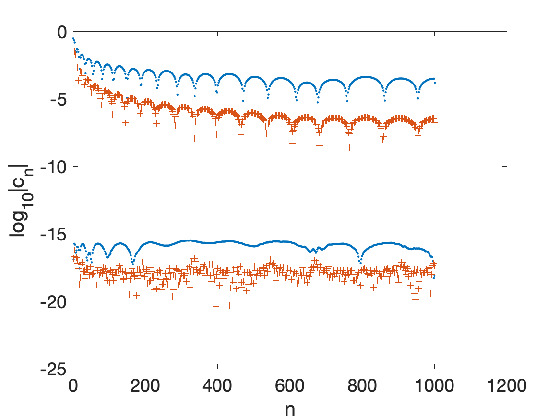}
 \includegraphics[width=0.32\textwidth]{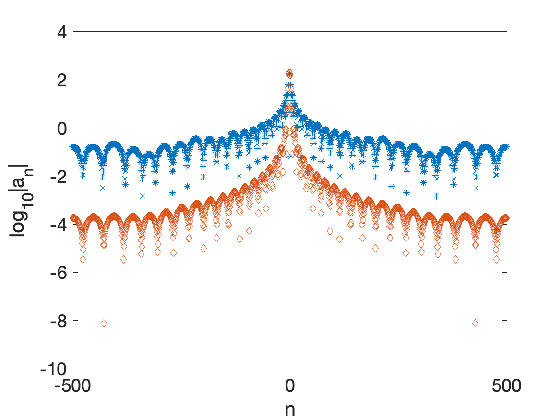}
 \includegraphics[width=0.32\textwidth]{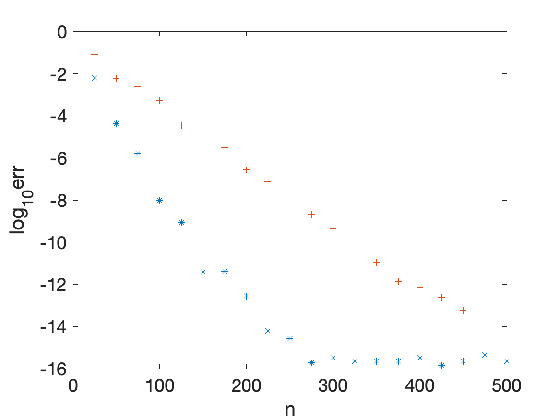}
 \caption{The spectral coefficients for the functions (\ref{sin12}) 
 for \( N_{\mathcal{F}}=1000 \), 
 on the left the Chebyshev coefficients in the infinite domain \( 
 |x|>1 \), in the middle the coefficients for the global approach 
 \cite{weideman}, in blue for \( f_{1} \), in red for \( f_{2} \);
 the numerical error for the computation of the Hilbert 
 transform for the functions (\ref{sin12}) with a contour deformation 
 approach are shown on the right.}
 \label{sinfig}
\end{figure}

As stated the multi-domain spectral approach is intended for 
functions piece-wise analytic on the whole real line or with rapid 
decrease towards infinity. It has been shown for various examples that it works as 
intended in such cases. In the case of an oscillatory behavior at 
infinity as for examples 3 and 4 of \cite{weideman}, 
\begin{equation}
	f_{1}(y) = \frac{\sin y}{1+y^{2}}, \quad f_{2}(y) = \frac{\sin 
	y}{1+y^{4}}
	\label{sin12},
\end{equation}
a spectral approximation is not ideal. We show the spectral 
coefficients both for the multi-domain approach in the infinite 
domain on the left and for the global approach \cite{weideman} in the 
middle of Fig.~\ref{sinfig}. The algebraic decay of the coefficients 
in both cases can be seen. 

The Hilbert transform for both functions is given by 
\begin{equation}
\begin{array}{c}
\mathcal{H}[f_{1}](x)=\frac{\cos(x)-\exp(-1)}{1+x^2},\\
~\\
	\mathcal{H}[f_{2}](x)=\frac{\cos(x)-\exp(-1/\sqrt{2})(\cos(1/\sqrt{2})+
	\sin(1/\sqrt{2})x^2)}{1+x^4}
\end{array}
	\label{hsin12}.
\end{equation}
The error for \( N_{\mathcal{F}}=1000 \) for the global 
approach is of the order of \( 10^{-3} \) for \( f_{1} \) and of the 
order of \( 10^{-7} \) for \( f_{2} \). To reach spectral 
convergence in such a case, deformation techniques in the complex 
plane appear to be necessary. Since the focus of this paper is 
on integration over
the real axis (in order to be able to deal with functions for which 
the localization of singularities in the complex plane is not known), 
this is in principle beyond the scope of the current paper. But we 
add this example for completeness and to show how the techniques can 
be efficiently extended in this way. For the examples (\ref{sin12}) 
we know that the singularities are, besides the obvious one on the 
real axis due to the Cauchy kernel, on the unit circle. Thus we can 
deform the integration contour from the real axis to \( y = 
e^{i\alpha }t+i\beta \), where \( \alpha,\beta \) are real constants 
and where \( t\in\mathbb{R} \). For an optimal choice of the deformed 
contours, \emph{steepest descent} techniques would have to be applied 
as in \cite{Tro,OT} in this context. 

Instead of integrating the sine function, we consider integrals of \( 
\exp(\pm iy )\) (or simply the imaginary part of the result for one 
of them) and choose for the first example in (\ref{sin12}) 
\( \alpha=\pm \pi/4 
\), \( \beta=\pm 0.5 \) and for the second \( \alpha=\pm\pi/8 \) and \( 
\beta=\pm0.2 \). The signs are always chosen in a way that the 
integrand is exponentially decreasing towards infinity on the 
considered interval. In this 
way we have mapped the problem to the case treated in 
Fig.~\ref{abserr}. We use the same parameters as there. Note that the 
terms proportional to \( \cos(x) \) in (\ref{hsin12}) are the contribution of 
the residue of the Cauchy kernel on the real axis which is thus taken 
care off analytically. On the deformed contours, the integrands are 
regular and no regularization as on the real axis is needed.
The numerical 
error for both examples can be seen on the right of 
Fig.~\ref{sinfig}. As expected spectral convergence is reached. 

\section{Solitary waves for generalized Benjamin-Ono equations}
Benjamin-Ono equations (\ref{BO}) appear for $m=2$ in applications for instance 
in the modelisation of two-layer fluids, see \cite{sautBO} and 
references therein. The case $m=2$ is in addition completely 
integrable. For $m> 2$, the solutions to initial value problems 
with smooth localized initial data of sufficiently large $L^{2}$ norm 
can have a blow-up in finite time and are thus mathematically 
interesting, see \cite{RWY} for a recent numerical study. 
As an application of the multi-domain spectral approach presented in 
the previous sections, we want to construct the solitary waves 
given numerically in \cite{RWY} 

A solitary wave is a traveling wave solution of (\ref{BO}) vanishing 
at infinity, i.e., a solution of the form $u(x,t)=Q_{c}(x-ct)$ where 
$c>0$ is a constant. Equation (\ref{BO}) implies for $Q$ 
the equation
\begin{equation}
	-cQ_{c}(\xi)-HQ_{c}'(\xi)+\frac{1}{m}Q_{c}^{m}(\xi) = 0,
	\label{Q}
\end{equation}
where we have integrated the equation resulting from (\ref{BO}) once 
using the vanishing of $Q$ at infinity; we have put $\xi=x-ct$ to 
stress that (\ref{Q}) is a nonlinear and nonlocal equation in one 
variable only. In addition we have the scaling invariance
$Q_{c}(\xi) = cQ(c\xi)$, where we have put $Q(\xi):=Q_{1}(\xi) $. 
Thus it is sufficient to consider the case $c=1$. The soliton in 
the integrable case $m=2$ is explicitly known,
\begin{equation}
	Q(\xi) = \frac{4}{1+\xi^{2}}, 
	\label{Qm2}
\end{equation}
the Lorentz profile we discussed in section 2 for the Hilbert transform. 

To numerically construct the solitary waves for $m>2$, we use the 
same approach as in section 2 with the two intervals $\xi\in[-1,1]$ and 
$1/\xi\in[-1,1]$, and the same for the computation of the Hilbert 
transform. For simplicity we use the same number $N_{1,\infty}$ of 
collocation points in both intervals for the integration in the 
variable $y$, but sample also $\xi$ on these points. Note that the 
Hilbert transform in the infinite interval is computed for the 
function $yQ(y)$. The derivative in (\ref{Q}) is approximated as 
before in terms of Chebyshev differentiation matrices. Thus we 
approximate (\ref{Q}) for $c=1$ by the discrete nonlinear equation system 
\begin{equation}
	\mathbf{F}(\mathbf{Q}):=-\mathbf{Q}-D\mathbb{H}\mathbf{Q}+\frac{1}{m}\mathbf{Q}^{m} = 0,
	\label{FQ}
\end{equation}
where $\mathbf{Q}$ is the vector with components of $Q(\xi_{n})$ with 
$\xi_{n}$, $n=0,\ldots,2N+1$ being the Chebyshev collocation points 
in the two intervals, 
where $\mathbb{H}$ is the matrix corresponding to the Hilbert 
transform, and where $D$ is the Chebyshev differentiation matrix as before. 
Thus we have to solve a nonlinear equation system which we do by a 
standard Newton iteration,
\begin{equation}
	\mathbf{Q}^{(k+1)} = \mathbf{Q}^{(k)} - 
	\mbox{Jac}^{-1}\mathbf{F}(\mathbf{Q}^{(k)})
	\label{newton},
\end{equation}
where $\mathbf{Q}^{(k)}$ is the $k$th iterate, $k=0,1,2,\ldots$, and 
where $\mbox{Jac}$ is the Jacobian of $\mathbf{F}(\mathbf{Q})$ with 
respect to $\mathbf{Q}$.

Note that the Jacobian has a kernel and thus cannot be inverted 
directly. This is partially due to a derivative appearing in the 
linear part of (\ref{Q}). To address this we require that $Q$ is 
continuous for $x=\pm1$. Both conditions are implemented with Lanczos' 
$\tau$-method \cite{tau}. In addition, equation (\ref{BO}) is translation 
invariant, if  $Q_{c}(\xi)$ is a solution so is $Q_{c}(\xi-\xi_{0})$ for constant $\xi_{0}$. To fix $\xi_{0}$, we require that 
$Q'(0)=0$. This we implement numerically with a $\tau$-method, to this end $N_{1,\infty}$ even in this section to make sure that $\xi=0$ is a 
collocation point).

As the zeroth iterate we use in all cases 
$A/(1+x^{2})$, $A>0$. The iteration is stopped once the $L^{\infty}$ 
norm of $\mathbf{F}$ is smaller than some threshold, typically 
$10^{-10}$. For $m>2$ we use some relaxation to stabilize the 
iteration, i.e., in each step of the iteration the new  iterate is 
formed by $\mu \mathbf{Q}^{(k+1)}+(1-\mu)\mathbf{Q}^{(k)}$ where 
$0<\mu\leq 1$.  

\begin{figure}[htb!]
 \includegraphics[width=0.49\textwidth]{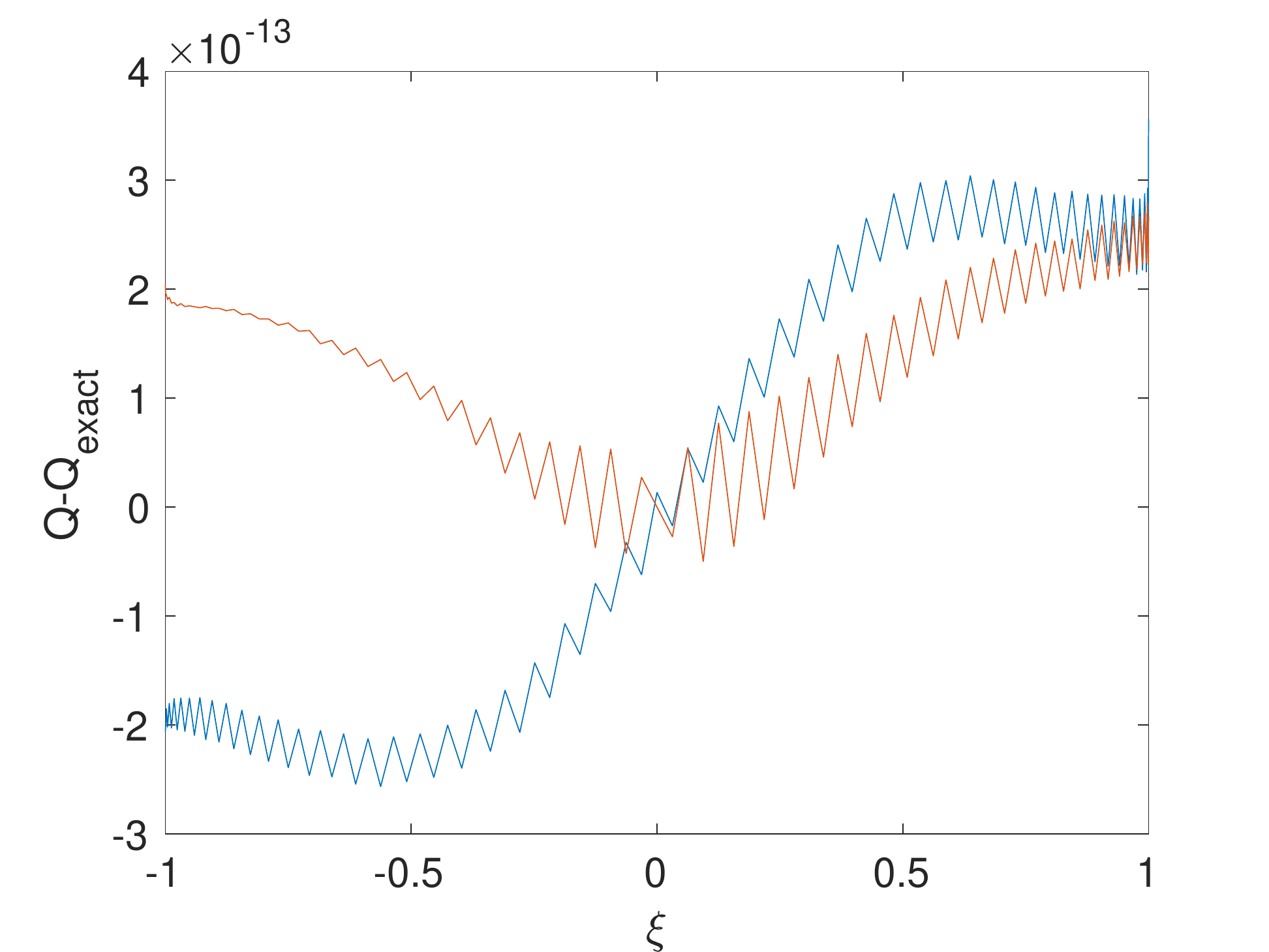}
 \includegraphics[width=0.49\textwidth]{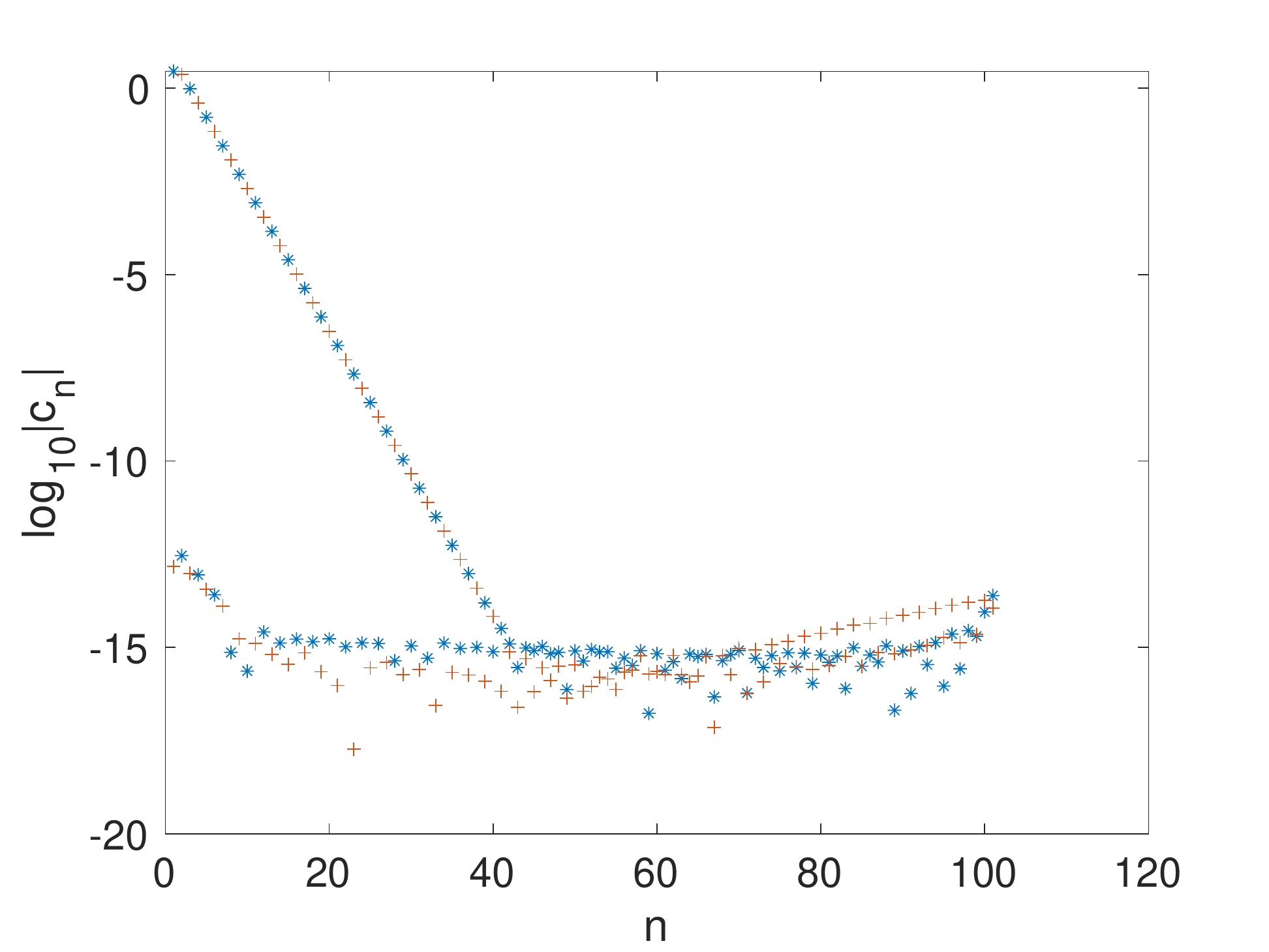}
  \caption{On the left the difference of the exact and the numerical 
  solution for the soliton for $m=2$, on the right the Chebyshev 
  coefficients of the solution, in blue for the finite domain, in red 
  for the infinite domain. }
 \label{M2fig}
\end{figure}

We first test the known solution for $m=2$ with $N_{1,\infty}=100$. For $A=3,5$ in the 
initial iterate, the Newton 
iteration converges after 5 iterations with a residual smaller than 
$10^{-10}$. The difference between numerical and exact solution can 
be seen on the left of Fig.~\ref{M2fig} to be of the order $10^{-13}$ 
in both domains. The Chebyshev coefficients of the solution in both 
domains on the right of the same figure indicate that $N\sim50$ 
collocation points are enough as in section 2 to reach maximal 
precision. 

\begin{figure}[htb!]
 \includegraphics[width=0.7\textwidth]{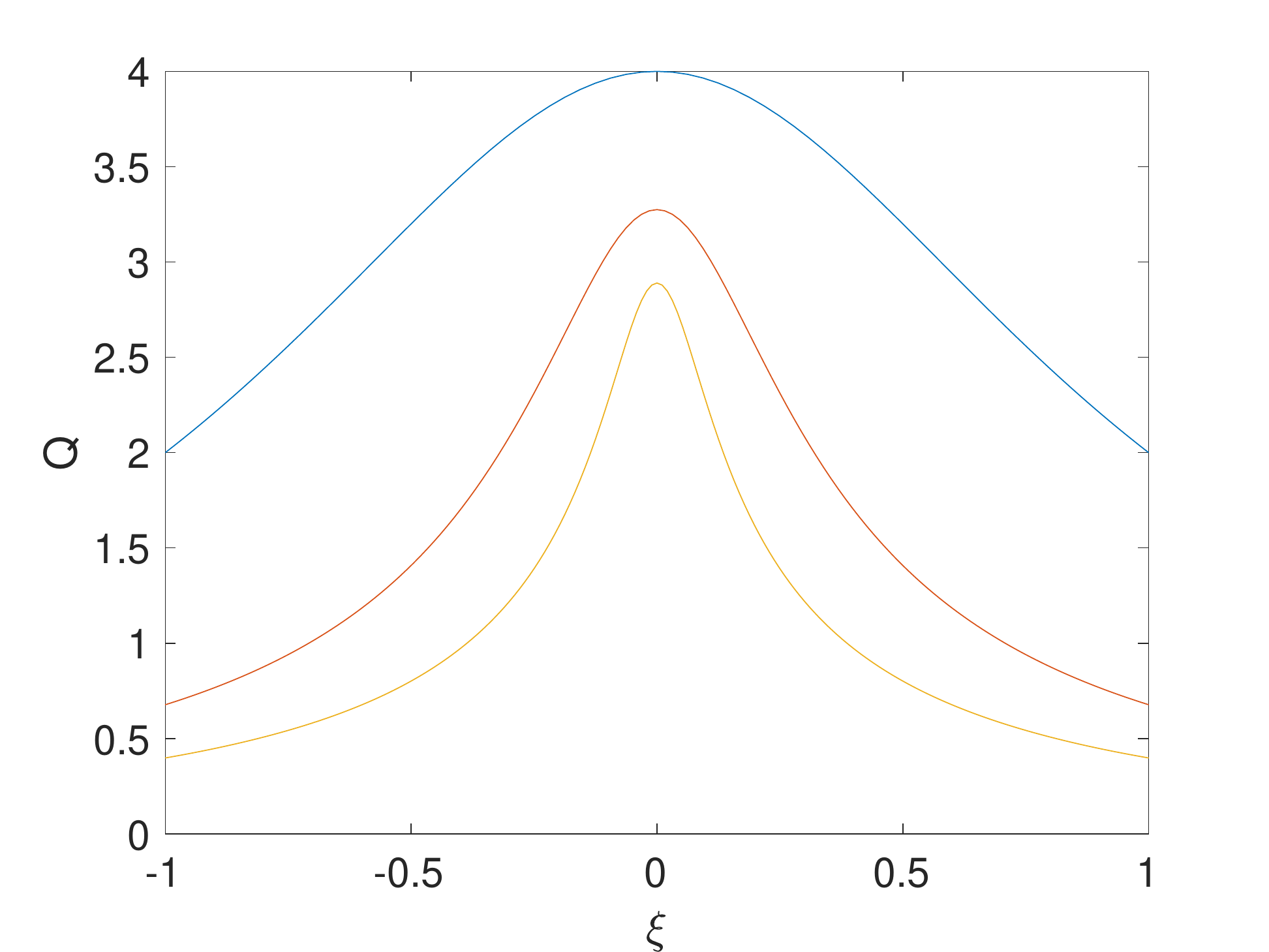}
  \caption{The solitary waves of (\ref{BO}) for $m=2,3,4$ (from top 
  to bottom). }
 \label{3Mfig}
\end{figure}

For higher values of $m$, we use $N_{1,\infty}=300$ collocation points in each 
domain. The solutions can be seen in Fig.~\ref{3Mfig} in the finite 
domain and are in accordance with \cite{RWY}. The higher the 
nonlinearity, the more the solitary waves are compressed. 

\begin{figure}[htb!]
 \includegraphics[width=0.49\textwidth]{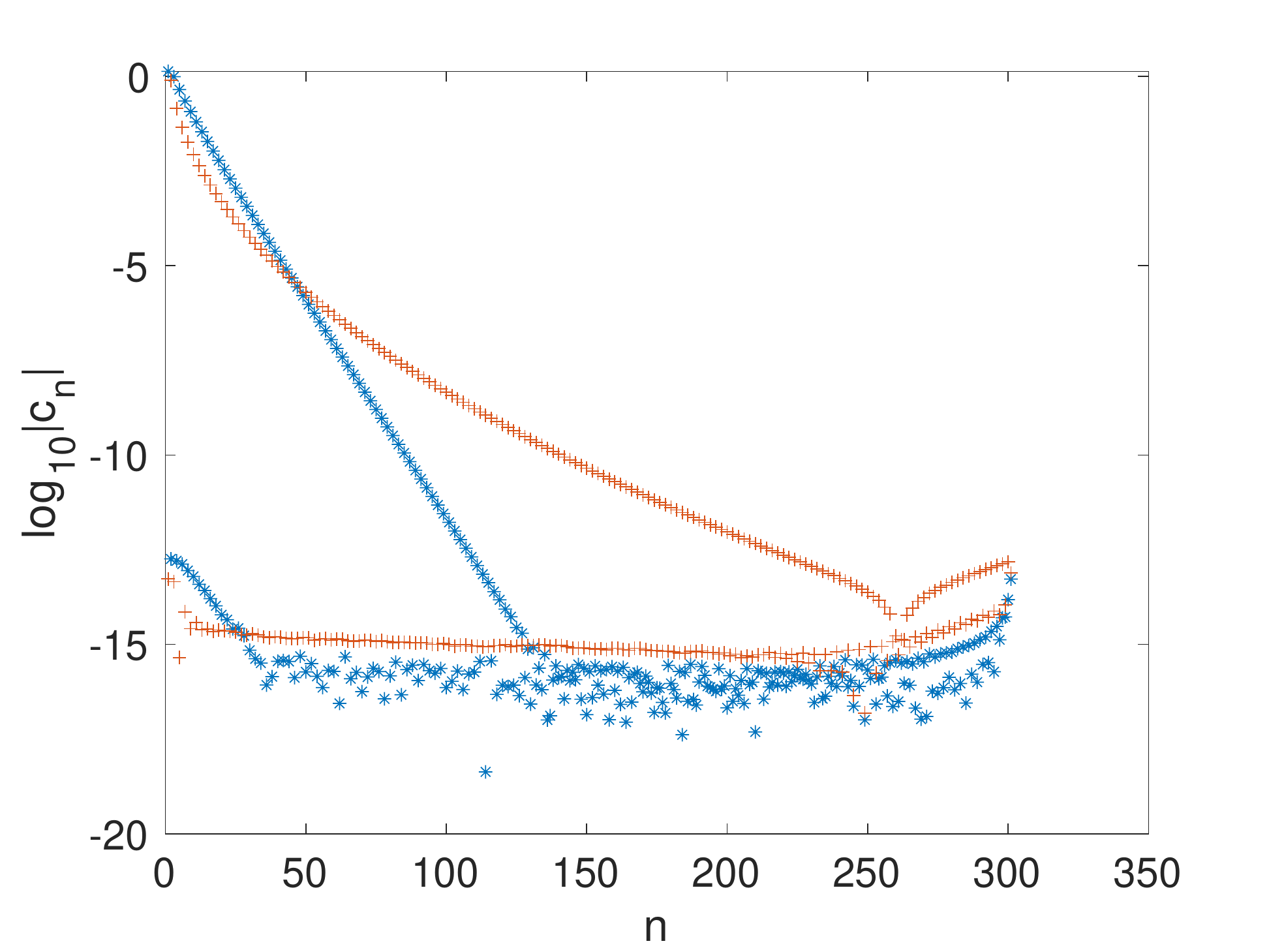}
 \includegraphics[width=0.49\textwidth]{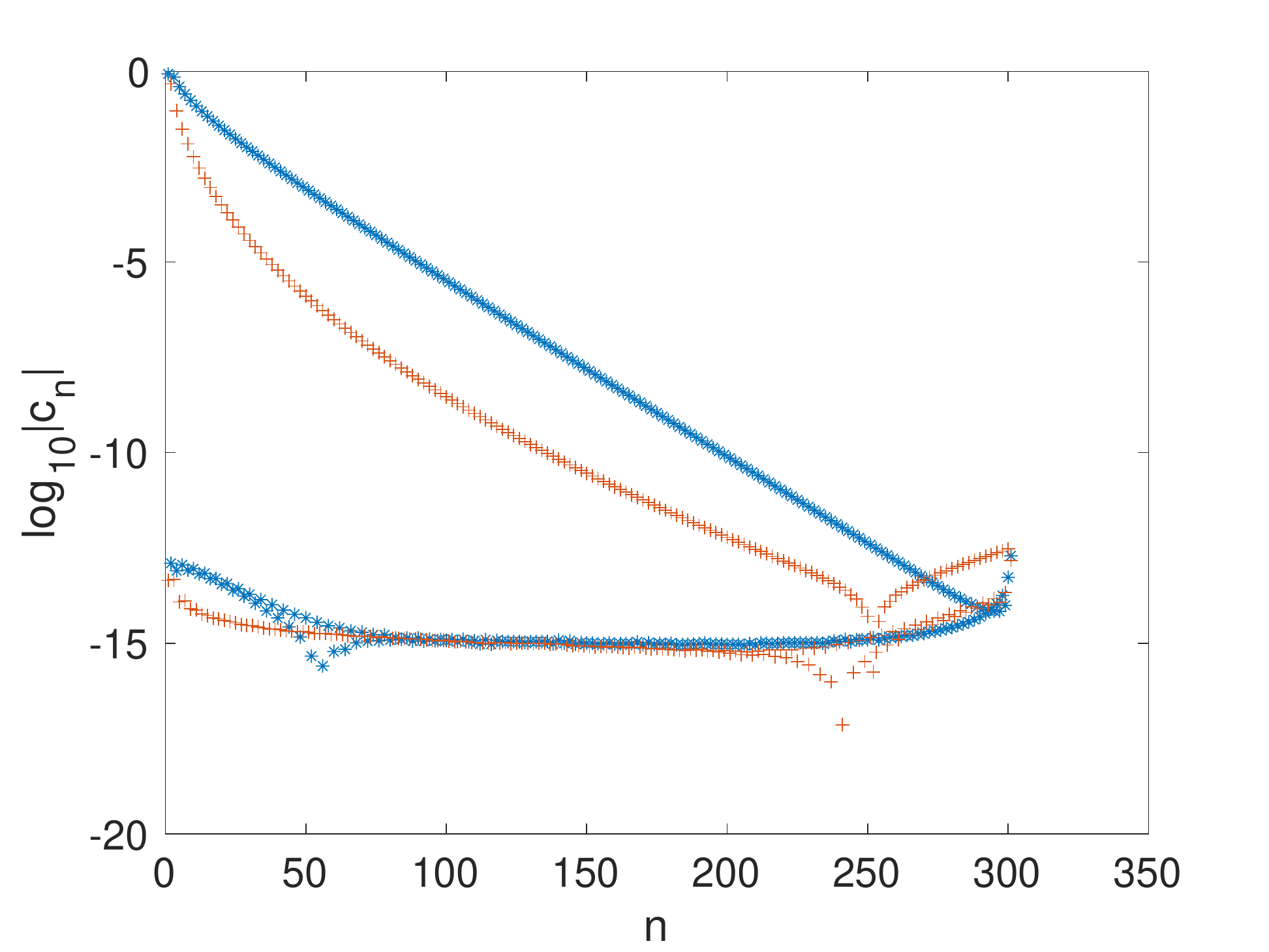}
  \caption{The Chebyshev 
  coefficients of the solitary waves for $m=3$ on the left and for 
  $m=4$ on the right, in blue for the finite domain, in red 
  for the infinite domain. }
 \label{2Mcoeff}
\end{figure}

The solutions for larger $m$ also require more numerical resolution, 
i.e., higher values of $N$ in each domain as can be seen in 
Fig.~\ref{2Mcoeff}. But for $N_{1,\infty}=300$, the coefficients decrease 
in all cases to the order of the rounding error. 

Note that the solitary waves have the symmetry $Q(-\xi)=Q(\xi)$. This 
is the reason why all odd Chebyshev coefficients in 
Fig.~\ref{2Mcoeff} vanish. To optimize resources, one could have 
worked on $\mathbb{R}^{+}$ only. But since we are in the future 
interested in studying the dynamics of perturbations of the solitary 
waves, i.e., use $Q$ plus perturbations as initial data for the 
generalized BO equation (\ref{BO}) as in \cite{RWY}, this is not 
convenient since the BO solution will not stay symmetric in general.

\section{Outlook}
In this paper we have presented a multi-domain spectral approach for 
the Hilbert transform on the real line. We have shown that it 
provides a comparable performance to Weideman's global approach 
\cite{weideman} for functions analytic on the whole real line. At 
various examples we have discussed that the global approach 
\cite{weideman} based on an expansion in terms of rational functions 
is more efficient for functions with an algebraic decrease towards 
infinity, but that this can be slightly different for rapidly 
decreasing functions. In all cases the same order of magnitude of 
collocation points is needed to achieve the same accuracy. The FFT 
based approach \cite{weideman} has a lower complexity and is thus the 
method of choice in such cases. The multi-domain approach is intended 
for piecewise analytic functions where it provides spectral 
accuracy when a global approach is of finite order and 
may exhibit Gibbs phenomena. This was illustrated at various examples.

One application of the multi-domain approach will be to study zones 
of rapid modulated oscillations called \emph{dispersive shock waves} 
(see for instance \cite{GK} for a review with many references) which appear in the solutions of nonlinear dispersive PDEs as the 
Benjamin-Ono equation (\ref{BO}). A multi-domain approach allows 
a special allocation of resolution where it is most needed, i.e., where the 
oscillations are. This is in particular interesting if one wants to 
study discontinuous initial data as in the case of the 
Gurevitch-Pitaevski \cite{GP} problem for the Korteweg-de Vries equation. 
Generalized BO equations (\ref{BO}) for sufficiently 
large \( p \) can have solutions to initial value problems with 
smooth initial data which blow up in finite time, i.e., where the \( 
L^{\infty}  \) norm diverges, see \cite{RWY} for a numerical study. The multi-domain approach will allow to 
study numerically such a blow-up with a combination of methods of 
\cite{BK} and a \emph{dynamical rescaling} as for the generalized 
Korteweg-de Vries equations in \cite{KP}. This will be the subject of 
further research. 

The multi-domain approach presented in the present paper was mainly 
developed for the real axis. However, as the example with an oscillatory 
singularity shows, it is straight forward to generalize this to 
arbitrary piecewise smooth contours in the complex plane. Each of 
the smooth arcs of such a contour (or parts of it) can be mapped to the interval \( 
[-1,1] \) where the same techniques as here can be applied to compute 
a Cauchy integral. The approach is also not limited to Cauchy 
integrals. In recent years, there has been an increasing interest in 
fractional derivatives, for instance in the context of PDEs with 
nonlocal dispersion, see e.g., \cite{KS,KSM,BLF}. The extent to which a multi-domain approach can be used to efficiently 
compute fractional derivatives will be studied in a separate work.

\end{document}